\documentclass[preprint,12pt]{elsarticle}
\usepackage{amssymb,amsmath,graphicx}

\journal{ }

\newtheorem{theo}{Theorem}
\newtheorem{defn}{Definition}

\newcommand{\ds}{\displaystyle}
\newcommand{\bc}{\begin{equation*}\begin{array}{l}}
\newcommand{\ec}{\end{array}\end{equation*}}

\begin{document}

\begin{frontmatter}

\title
{A Family of Multistep Methods with Zero Phase-Lag and Derivatives for the Numerical Integration of Oscillatory ODEs}

\author[UoP]{Z. A. Anastassi}
\ead{zackanas@uop.gr}

\author[UoP]{D. S. Vlachos}
\ead{dvlachos@uop.gr}

\author[UoP]{T. E. Simos\fnref{simos}}
\fntext[simos]{Highly Cited Researcher, Active Member of the European Academy of Sciences and Arts, Address: Dr. T.E. Simos, 26 Menelaou Street, Amfithea - Paleon Faliron, GR-175 64 Athens, GREECE, Tel: 0030 210 94 20 091}
\ead{tsimos.conf@gmail.com, tsimos@mail.ariadne-t.gr}

\address[UoP]{Laboratory of Computer Sciences,\\
Department of Computer Science and Technology,\\
Faculty of Sciences and Technology, University of Peloponnese\\
GR-22 100 Tripolis, GREECE}

\begin{abstract}
\noindent In this paper we develop a family of three 8-step methods, optimized for the numerical integration of oscillatory ordinary differential equations. We have nullified the phase-lag of the methods and the first $r$ derivatives, where $r=\{1,2,3\}$. We show that with this new technique, the method gains efficiency with each derivative of the phase-lag nullified. This is the case for the integration of both the Schr\"odinger equation and the N-body problem. A local truncation error analysis is performed, which, for the case of the Schr\"odinger equation, also shows the connection of the error and the energy, revealing the importance of the zero phase-lag derivatives. Also the stability analysis shows that the methods with more derivatives vanished, have a bigger interval of periodicity.
\end{abstract}


\end{frontmatter}

\section{Introduction}
The numerical integration of systems of ordinary differential equations with oscillatory solutions has been the subject of research during the past decades. This type of ODEs is often met in real problems, like the N-body problem and the Schr\"odinger equation.

There are some special techniques for optimizing numerical methods. Trigonometrical fitting and phase-fitting are some of them, producing methods with variable coefficients, which depend on $v=\omega h$, where $\omega$ is the dominant frequency of the problem and $h$ is the step length of integration.

For example Raptis and Allison have developed a two-step exponentially-fitted method of order four in
\cite{rallison} and Kalogiratou and Simos have constructed a two-step P-stable exponentially-fitted method of
order four in \cite{kalogiratou}. Also Panopoulos, Anastassi and Simos have constructed two optimized eight-step methods with high or infinite order of phase-lag in \cite{panopoulos_match}.

Some other notable multistep methods for the numerical solution of oscillating IVPs have been developed by Chawla and
Rao in \cite{chawla}, who produced a three-stage, two-Step P-stable method with minimal phase-lag and order six and by
Henrici in \cite{henrici}, who produced a four-step symmetric method of order six. Also some recent research work in numerical methods can be found in \cite{anastassi1}, \cite{anastassi2}, \cite{anastassi3}, \cite{Meyer}, \cite{ref_12}, \cite{ref_14}, \cite{vanden_jnaiam}, \cite{cash_jnaiam1}, \cite{cash_jnaiam2}, \cite{iavernaro_jnaiam1}, \cite{simos_cole1}, \cite{simos_cole2} and \cite{Psihoyios_cole}.

Trigonometrically fitted methods of high trigonometric order are well known for their high efficiency in the integration of the Schr\"odinger equation, especially when using a high value of energy. However higher trigonometric order is not rendering them more efficient for all types of oscillatory problems. On the other hand, phase-lag does not give us the opportunity to provide such methods, that for example perform well when integrating the Schr\"odinger equation for high values of energy.

In this paper we present a methodology for optimizing numerical methods, through the use of phase-lag and its derivatives with respect to $v$. More specifically, given a classical (i.e. with constant coefficients) numerical method, we can provide a family of optimized methods, each of which has zero $\{PL\}$ or zero $\{PL$ and $PL'\}$ or zero $\{PL$, $PL'$ and $PL''\}$ etc.

With this new technique we provide methods that perform well during the integration of the Schr\"odinger equation for high values of energy, but also that perform well on other real problems with oscillatory solution, like the N-body problem.

\section{Phase-lag and stability analysis of symmetric multistep methods}

For the numerical solution of the initial value problem

\begin{equation}
\label{ivp_definition}
    y'' = f(x,y)
\end{equation}

\noindent multistep methods of the form

\begin{equation}
\label{multistep_definition}
    \sum\limits_{i=0}^{m}{a_{i}y_{n+i}} = h^{2}\sum\limits_{i=0}^{m}{b_{i}f(x_{n+i},y_{n+i})}
\end{equation}

with $m$ steps can be used over the equally spaced intervals $\left\{x_{i}\right\}^{m}_{i=0} \in [a,b]$ and
$h=|x_{i+1}-x_{i}|$, \, $i=0(1)m-1$.

If the method is symmetric then $a_i=a_{m-i}$ and $b_i=b_{m-i}$, \, $i=0(1)\lfloor \frac{m}{2} \rfloor$.

Method \eqref{multistep_definition} is associated with the operator

\begin{eqnarray}
\label{exp_operator} L(x) = \sum\limits_{i=0}^{m}{a_{i}u(x+ih)} - h^{2}\sum\limits_{i=0}^{m}{b_{i}u''(x+ih)}
\end{eqnarray}

\noindent where $u \in C^2$.

\begin{defn}
\label{defn_exp1} \emph{The multistep method \eqref{exp_operator} is called algebraic of order $p$ if the associated
linear operator $L$ vanishes for any linear combination of the linearly independent functions $\,1,\, x,\,x^2,\, \ldots
,\, x^{p+1}$.}
\end{defn}

When a symmetric $2k$-step method, that is for $i=-k(1)k$, is applied to the scalar test equation

\begin{equation}
\label{stab_eq} y''=-\omega^2 y
\end{equation}

a difference equation of the form
\begin{eqnarray}
\label{phl_multi_de}
\nonumber A_{k} (v)y_{n + k} + ... + A_{1} (v)y_{n + 1} + A_{0} (v)y_{n}\\
+ A_{1}(v)y_{n - 1} + ... + A_{k} (v)y_{n - k} = 0
\end{eqnarray}

\noindent is obtained, where $v = \omega h$, $h$ is the step length and $A_{0} (v)$, $A_{1} (v),\ldots$, $ A_{k} (v)$
are polynomials of $v$.

The characteristic equation associated with \eqref{phl_multi_de} is

\begin{eqnarray}
\label{phl_multi_ce}
A_{k} (v)s^{k} + ... + A_{1} (v)s + A_{0} (v) + A_{1} (v)s^{ - 1} + ... + A_{k} (v)s^{ - k} = 0
\end{eqnarray}

\begin{theo}
\emph{\cite{royal}} The symmetric $2k$-step method with characteristic equation given by \eqref{phl_multi_ce} has
phase-lag order $q$ and phase-lag constant $c$ given by

\begin{equation}
\label{phl_multi_defn} - c v ^{q + 2} + O(v^{q + 4}) = {\frac{{2A_{k} (v)\cos (k v ) + ... + 2A_{j} (v)\cos (j v ) + ...
+ A_{0} (v)}}{{2k^{2}A_{k} (v) + ... + 2j^{2}A_{j} (v) + ... + 2A_{1} (v)}}}
\end{equation}
\end{theo}

The formula proposed from the above theorem gives us a direct method to calculate the phase-lag of any symmetric $2k$-
step method.

The characteristic equation has $m$ characteristic roots $\lambda_{i}, \; i=0(1)m-1$.

\begin{defn} \emph{
\cite{la_wa} If the characteristic roots satisfy the conditions $|\lambda_{i}|\leqslant 1, \; i=0(1)m-1$ for all $s=\theta h$, then we say that the method is} unconditionally stable.
\end{defn}

\begin{defn} \label{sta_int} \emph{
\cite{la_wa} If the characteristic roots satisfy the conditions $\lambda_1=e^{I\,\phi(s)}$, $\lambda_2=e^{-I\,\phi(s)}$ $|\lambda_{i}|\leqslant 1, \; i=3(1)m-1$ for all $s<s_{0}$, where $s=\theta h$ and $\phi(s)$ is a real function of $s$, then we say that the method has interval of periodicity $(0,s_{0}^2)$.}
\end{defn}

\begin{defn} \label{p_stability} \emph{
\cite{la_wa} Method \eqref{multistep_definition} is called P-stable if its \emph{interval of periodicity} is $(0,\infty)$.}
\end{defn}

\section{Construction of the new optimized multistep methods}
\label{Construction}

We consider the multistep symmetric method of Quinlan-Tremaine \cite{qt8}, with eight steps and eighth algebraic order:

\begin{equation}
\begin{array}{c} \label{table_qt8}
y_{{4}} = -y_{{-4}} -a_{{3}}(y_{{3}}+y_{{-3}}) -a_{{2}}(y_{{2}}+y_{{-2}}) -a_{{1}}(y_{{1}}+y_{{-1}})\\
+{h}^{2}\left(b_{{3}}(f_{{3}}+f_{{-3}}) +b_{{2}}(f_{{2}}+f_{{-2}}) +b_{{1}}(f_{{1}}+f_{{-1}}) +b_{{0}}f_{{0}}\right)
\end{array}
\end{equation}

\noindent where

\begin{equation}
\begin{array}{l}
a_{3}=-2, \qquad a_{2}=2, \qquad a_{1}=-1, \vspace{5pt}\\
\ds b_{3}=\frac{17671}{12096}, \qquad b_{2}=-\frac{23622}{12096}, \qquad
b_{1}=\frac{61449}{12096}, \qquad b_{0}=-\frac{50516}{12096},\vspace{5pt}\\
y_{i} = y(x+ih) \mbox{ and } f_{i} = f(x+ih,y(x+ih))
\end{array}
\end{equation}

We also consider the optimized method, that is based on the above one, with zero phase-lag constructed by Panopoulos, Anastassi and Simos in \cite{panopoulos_match}. The coefficients are given below:

\begin{equation}
\label{meth_PL_0}
\begin{array}{l}
\ds b_{{0}}=-20\,b_{{3}}+{\frac {601}{24}}, \qquad b_{{2}}=-6\,b_{{3}}+{ \frac {109}{16}}, \qquad b_{{1}}=15\,b_{{3}}-{\frac{101}{6}}\\
\ds \nonumber b_{3} = {\frac {1}{96}}\,{\frac {C}{D}}, \qquad \mbox{where}\\
\nonumber C=-192\, \left( \cos \left( v \right)  \right) ^{4}+192\, \left( \cos \left( v \right) \right) ^{3}+ \left(
96-327\,{v}^{2} \right)  \left( \cos \left( v \right) \right) ^{2}\\
\nonumber + \left( - 120+404\,{v}^{2} \right) \cos \left( v \right) -137\,{v}^{2}+24\\
D={{v}^{ 2} \left( \cos \left( v \right) -1 \right) ^{3}},
\end{array}
\end{equation}

\noindent where $v=\omega h$ and the $a_{i}$ coefficients remain the same. The Taylor series expansions of the coefficients are:

\bc
\ds b_{{0}}=-{\frac {12629}{3024}}+{\frac {45767}{36288}}\,{v}^{2}-{\frac{164627}{2395008}}\,{v}^{4}+{\frac {520367}{792529920}}\,{v}^{6}\vspace{5pt}\\
\ds -{\frac {76873}{4483454976}}\,{v}^{8}+{\frac {9190171}{160059342643200}}\,{v}^{10}+{\frac {6662921}{1703031405723648}}\,{v}^{12}+\ldots\\
\\
\ds b_{{1}}={\frac {20483}{4032}}-{\frac {45767}{48384}}\,{v}^{2}+{\frac {164627}{3193344}}\,{v}^{4}-{\frac {520367}{1056706560}}\,{v}^{6}\vspace{5pt}\\
\ds +{\frac {76873}{5977939968}}\,{v}^{8}-{\frac {9190171}{213412456857600}}\,{v}^{10}-{\frac {6662921}{2270708540964864}}\,{v}^{12}-\ldots\\
\\
\ds b_{{2}}=-{\frac {3937}{2016}}+{\frac {45767}{120960}}\,{v}^{2}-{\frac{164627}{7983360}}\,{v}^{4}+{\frac {520367}{2641766400}}\,{v}^{6}\vspace{5pt}\\
\ds -{\frac {76873}{14944849920}}\,{v}^{8}+{\frac {9190171}{533531142144000}}\,{v}^{10}+{\frac {6662921}{5676771352412160}}\,{v}^{12}+\ldots\\
\\
\ds b_{{3}}={\frac {17671}{12096}}-{\frac {45767}{725760}}\,{v}^{2}+{\frac {164627}{47900160}}\,{v}^{4}-{\frac {520367}{15850598400}}\,{v}^{6}\vspace{5pt}\\
\ds +{\frac {76873}{89669099520}}\,{v}^{8}-{\frac {9190171}{3201186852864000}}\,{v}^{10}-{\frac {6662921}{34060628114472960}}\,{v}^{12}-\ldots
\ec

We want to produce three new methods that, apart from zero phase-lag, will also have zero $r$ derivatives of the phase-lag, where $r=\{1,2,3\}$. In particular the three new methods must satisfy these equations:

\begin{itemize}
\item First method:  $\{PL=0,PL'=0\}$
\item Second method: $\{PL=0,PL'=0,PL''=0\}$
\item Third method:  $\{PL=0,PL'=0,PL''=0,PL'''=0\}$
\end{itemize}

Since we have four free coefficients $b_i$, $i=\{0,1,2,3\}$ ($a_{i}$ remain the same), the rest of the coefficients for each method will be determined by the algebraic conditions.

\subsection{First optimized method with zero $PL$ and $PL'$}

The first method must satisfy the conditions $\{PL=0,PL'=0\}$, thus we need two coefficients to be determined by the maximum algebraic order.

We use formula \eqref{phl_multi_defn} to compute the phase-lag and then its first derivative in respect to $v$:

\noindent where $v=\omega\,h$, $\omega$ is the frequency and $h$ is the step length used.

\bc
PL = \big(96\, \left( \cos \left( v \right)  \right) ^{4}+ \left( -96+48 \,{v}^{2}b_{{3}} \right)  \left( \cos \left( v \right)  \right) ^{3}+\\
  \left( -48+24\,{v}^{2}b_{{2}} \right)  \left( \cos \left( v \right) \right) ^{2}
+ \left( 60+ \left( 125-144\,b_{{3}}-48\,b_{{2}} \right) {v}^{2} \right) \cos \left( v \right)\\
 -12+ \left( 24\,b_{{2}}-95+96\,b_{{3}} \right) {v}^{2}\big)\\
/\left({60+125\,{v}^{2}}\right)\\
\\
\ds PL' = \frac{1}{5}\, \Big( -4800\,v \left( \cos \left( v \right)  \right) ^{4}+ \big( 1152\,b_{{3}}v-9600\,\sin \left( v \right) {v}^{2}+4800\,v\\
\ds -4608\,\sin \left( v \right)  \big)  \left( \cos \left( v \right) \right) ^{3}+ \left( -3600\, \left( {\frac {12}{25}}+{v}^{2} \right) \left( -2+{v}^{2}b_{{3}} \right) \sin \left( v \right)\right.\\
\ds \left. + \left( 576\,b_{{2}}+2400 \right) v \right)  \left( \cos \left( v \right)  \right)^{2}+ \left( -1200\, \left( {\frac {12}{25}}+{v}^{2} \right)  \left( {v}^{2}b_{{2}}-2 \right) \sin \left( v \right)\right.\\
\ds \left. -3456\,v \left( b_{{3}}+1/3\,b_{{2}} \right)  \right) \cos \left( v \right)\\
\ds +3600\, \left( {\frac {12}{25}}+{v}^{2} \right)  \left( -{\frac {5}{12}}+ \left( -{\frac {125}{144}}+b_{{3}}+\frac{1}{3}\,b_{{2}} \right) {v}^{2} \right) \sin \left( v \right)\\
 +2304\,v \left( -{\frac {35}{48}}+b_{{3}}+1/4\,b_{{2}} \right)  \Big)  \left( 12+25\,{v}^{2} \right) ^{-2}
\ec

The four equations to be solved are:

$$
PL=0, \quad PL'=0, \quad b_0 = -{\frac {95}{6}}+16\,b_{{3}}+6\,b_{{2}}, \quad
b_1 = {\frac {125}{12}}-9\,b_{{3}}-4\,b_{{2}}
$$

\noindent and the coefficients are given below:

\begin{equation}
\label{coeff_meth_phl_deriv_1}
\begin{array}{l}
\ds b_0=-\frac{b_{0,num}}{12 D}, \quad b_1=\frac{b_{1,num}}{48 D}, \quad b_2=-\frac{b_{2,num}}{24 D}, \quad b_3=\frac{b_{3,num}}{48 D}\vspace{5pt},\\
\mbox{where } D = {\left(\left(\cos\left(v\right)\right)^{4}-2\,\left(\cos
\left(v\right)\right)^{3}+2\,\cos\left(v\right)-1\right){v
}^{3}}\;\;\; \mbox{and}
\end{array}
\end{equation}

\bc
b_{{0,num}}=-288\,\left(\cos\left(v\right)\right)^{
6}v+576\,\sin\left(v\right)\left(\cos\left(v\right)\right)
^{5}+192\,\left(\cos\left(v\right)\right)^{5}v\\
-192\,\sin
\left(v\right)\left(\cos\left(v\right)\right)^{4}+190\,
\left(\cos\left(v\right)\right)^{4}{v}^{3}+720\,\left(\cos
\left(v\right)\right)^{4}v\\
-120\,\left(\cos\left(v\right)
\right)^{3}v-672\,\sin\left(v\right)\left(\cos\left(v
\right)\right)^{3}+370\,\left(\cos\left(v\right)\right)^{3
}{v}^{3}\\
+168\,\sin\left(v\right)\left(\cos\left(v\right)
\right)^{2}-540\,\left(\cos\left(v\right)\right)^{2}v+145\,
\left(\cos\left(v\right)\right)^{2}{v}^{3}\\
+168\,\cos\left(v
\right)\sin\left(v\right)-70\,\cos\left(v\right){v}^{3}-72\,
\cos\left(v\right)v+108\,v\\
-48\,\sin\left(v\right) -35\,{v}^{3}
\ec

\bc
b_{{1,num}}=-768\,\left(\cos\left(v\right)\right)^{6
}v+1536\,\sin\left(v\right)\left(\cos\left(v\right)\right)
^{5}+192\,\left(\cos\left(v\right)\right)^{5}v\\
+500\,\left(
\cos\left(v\right)\right)^{4}{v}^{3}-192\,\sin\left(v\right)
\left(\cos\left(v\right)\right)^{4}+2400\,\left(\cos\left(
v\right)\right)^{4}v\\
+1000\,\left(\cos\left(v\right)\right)
^{3}{v}^{3}-2112\,\sin\left(v\right)\left(\cos\left(v\right)
\right)^{3}-1980\,\left(\cos\left(v\right)\right)^{2}v\\
+595\,
\left(\cos\left(v\right)\right)^{2}{v}^{3}+288\,\sin\left(v
\right)\left(\cos\left(v\right)\right)^{2}-100\,\cos\left(
v\right){v}^{3}\\
+648\,\cos\left(v\right)\sin\left(v\right)-
192\,\cos\left(v\right)v+348\,v-195\,{v}^{3}-168\,\sin\left(v
\right)
\ec

\bc
b_{{2,num}}=-96\,\left(\cos\left(v\right)\right)^{6
}v+192\,\sin\left(v\right)\left(\cos\left(v\right)\right)^
{5}-192\,\left(\cos\left(v\right)\right)^{5}v\\
+192\,\sin
\left(v\right)\left(\cos\left(v\right)\right)^{4}+624\,
\left(\cos\left(v\right)\right)^{4}v+216\,\left(\cos\left(
v\right)\right)^{3}v\\
-480\,\sin\left(v\right)\left(\cos
\left(v\right)\right)^{3}+250\,\left(\cos\left(v\right)
\right)^{3}{v}^{3}-612\,\left(\cos\left(v\right)\right)^{2}v\\
+215\,\left(\cos\left(v\right)\right)^{2}{v}^{3}-72\,\sin
\left(v\right)\left(\cos\left(v\right)\right)^{2}-70\,\cos
\left(v\right){v}^{3}\\
+216\,\cos\left(v\right)\sin\left(v
\right)-24\,\cos\left(v\right)v+84\,v-48\,\sin\left(v\right)
-35\,{v}^{3}
\ec

\bc
\ds b_{{3,num}}=\frac{1}{48}-192\,\left(\cos\left(v\right)\right)^{5
}v+192\,\sin\left(v\right)\left(\cos\left(v\right)\right)^
{4}+288\,\left(\cos\left(v\right)\right)^{4}v\\
+192\,\left(
\cos\left(v\right)\right)^{3}v-192\,\sin\left(v\right)
\left(\cos\left(v\right)\right)^{3}-96\,\sin\left(v\right)
\left(\cos\left(v\right)\right)^{2}\\
-324\,\left(\cos\left(v
\right)\right)^{2}v+125\,\left(\cos\left(v\right)\right)^{
2}{v}^{3}+120\,\cos\left(v\right)\sin\left(v\right)+60\,\cos
\left(v\right){v}^{3}\\
-24\,\sin\left(v\right) +36\,v-65\,{v}^{3}
\ec

The Taylor series expansions, used when $v \rightarrow 0$, are given below:

\bc
\ds b_{{0}}=-{\frac{12629}{3024}}+{\frac{45767}{18144}}\,{v}^{2}-{\frac{11483491}{23950080}}\,{v}^{4}+{\frac{112258001}{2615348736}}\,{v}^{6}\vspace{5pt}\\
\ds -{\frac{1703481341}{784604620800}}\,{v}^{8}+{\frac{5614773343}{80029671321600}}\,{v}^{10}-{\frac{10940565121}{6307523724902400}}\,{v}^{12}+...\\
\\
\ds b_{{1}}={\frac{20483}{4032}}-{\frac{45767}{24192}}\,{v}^{2}+{\frac{10476617}{31933440}}\,{v}^{4}-{\frac{45578707}{1585059840}}\,{v}^{6}\vspace{5pt}\\
\ds +{\frac{1514526707}{1046139494400}}\,{v}^{8}-{\frac{5016343559}{106706228428800}}\,{v}^{10}+{\frac{19742264573}{17466988776652800}}\,{v}^{12}-\\
\\
\ds b_{{2}}=-{\frac{3937}{2016}}+{\frac{45767}{60480}}\,{v}^{2}-{\frac{1491199}{15966720}}\,{v}^{4}+{\frac{321593093}{43589145600}}\,{v}^{6}\vspace{5pt}\\
\ds -{\frac{189532561}{523069747200}}\,{v}^{8}+{\frac{460150601}{38109367296000}}\,{v}^{10}-{\frac{28082396599}{113535427048243200}}\,{v}^{12}+...\\
\\
\ds b_{{3}}={\frac{17671}{12096}}-{\frac{45767}{362880}}\,{v}^{2}+{\frac{96865}{19160064}}\,{v}^{4}-{\frac{21971953}{261534873600}}\,{v}^{6}\vspace{5pt}\\
\ds +{\frac{82561}{448345497600}}\,{v}^{8}-{\frac{17608099}{123122571264000}}\,{v}^{10}-{\frac{1184824691}{75690284698828800}}\,{v}^{12}-...
\ec

\subsection{Second optimized method with zero $PL$, $PL'$ and $PL''$}

The second method must satisfy the conditions $\{PL=0, PL'=0, PL''=0\}$, thus we need one coefficient to be determined by the maximum algebraic order.

We use formula \eqref{phl_multi_defn} to compute the phase-lag and then its first and second derivative in respect to $v$:

\bc
\ds PL = \big(16\,\left(\cos\left(v\right)\right)^{4}+\left( 8\,b_{{3}}{v}^{2}-16\right)\left(\cos\left(v\right)\right)^{3}+
\left(4\,b_{{2}}{v}^{2}-8\right)\left(\cos\left(v\right)
\right)^{2}\\
+\left(10+\left(-6\,b_{{3}}+2\,b_{{1}}\right) {v}^{2}\right)\cos\left(v\right)-2+\left( -4\,b_{{2}}-2\,b_{{1}}-2\,b
_{{3}}+5\right){v}^{2}\big)\\
/\left(10+\left( 18\,b_{{3}}+8\,b_{{2}}+2\,b_{{1}}
\right){v}^{2}\right)
\ec

\bc
\ds PL' = \left(-16\,v\left(9\,b_{{3}}+4\,b_{{2}}+b_{{1}} \right)\left(\cos\left(v\right)\right)^{4}+\left(\left(-160+\left(-32\,
b_{{1}}-128\,b_{{2}}\right.\right.\right.\right.\\
\left.\left.\left.\left.-288\,b_{{3}}\right){v}^{2}\right)\sin\left(
v\right)+16\,\left(\frac{23}{2}\,b_{{3}}+4\,b_{{2}}+b_{{1}}\right)v
\right)\left(\cos\left(v\right)\right)^{3}+\right.\\
\left.\left(-12\,
\left(5+\left(9\,b_{{3}}+4\,b_{{2}}+b_{{1}}\right){v}^{2}
\right)\left(-2+b_{{3}}{v}^{2}\right)\sin\left(v\right)\right.\right.\\
\left.\left.+8\,\left(9\,b_{{3}}+\frac{13}{2}\,b_{{2}}+b_{{1}}\right)v\right)\left(
\cos\left(v\right)\right)^{2}+\left(-4\,\left(b_{{2}}{v}^{2}
-2\right)\left(5+\right.\right.\right.\\
\left.\left.\left.\left(9\,b_{{3}}+4\,b_{{2}}+b_{{1}}\right) {v}^{2}\right)\sin\left(v\right)-40\,v\left(3\,b_{{3}}+b_{{2}}
\right)\right)\cos\left(v\right)-\right.\\
\left.\left(5+\left(-3\,b_{{3}}
+b_{{1}}\right){v}^{2}\right)\left(5+\left( 9\,b_{{3}}+4\,b_{{2}}+b_{{1}}\right){v}^{2}\right)\sin\left(v\right)\right.\\
\left.-8\,v\left(
\frac{3}{2}\,b_{{2}}-b_{{3}}+b_{{1}}-{\frac{25}{8}}\right)\right)\\ /\left(5+\left(9\,b_{{3}}+4\,b_{{2}}+b_{{1}}\right){v}^{2}\right)^{2}
\ec

\bc
PL'' = {\frac{1}{729}}\,\left(\left(-10368\,\left( b_{{3}}+\frac{4}{9}\,b_{{2}}
+\frac{1}{9}\,b_{{1}}\right)^{2}{v}^{4}\right.\right.\\
\left.\left.+3888\,\left( \frac{1}{9}\,b_{{1}}-{\frac{
80}{27}}+b_{{3}}+\frac{4}{9}\,b_{{2}}\right)\left( b_{{3}}+\frac{4}{9}\,b_{{2}}+\frac{1}{9}
\,b_{{1}}\right) {v}^{2}-3200-720\,b_{{3}}\right.\right.\\
\left.\left.-320\,b_{{2}}-80\,b_{{1}}
\right)\left(\cos\left(v\right)\right)^{4}+\left(10368\,
\left(\frac{5}{9}+\left(b_{{3}}+\frac{4}{9}\,b_{{2}}+\frac{1}{9}\,b_{{1}}\right) {v}^{2}
\right)v\right.\right.\\
\left.\left.\left(b_{{3}}+\frac{4}{9}\,b_{{2}}+\frac{1}{9}\,b_{{1}}\right)\sin
\left(v\right)-2916\,\left(b_{{3}}+\frac{4}{9}\,b_{{2}}+\frac{1}{9}\,b_{{1}}
\right)^{2}b_{{3}}{v}^{6}\right.\right.\\
\left.\left.+2592\,\left( b_{{3}}+\frac{4}{9}\,b_{{2}}+\frac{1}{9}\,b_
{{1}}\right)\left(b_{{2}}+b_{{3}}+\frac{1}{4}\,b_{{1}}\right){v}^{4}+
\left(-768\,{b_{{2}}}^{2}\right.\right.\right.\\
\left.\left.\left.+\left( 2880-3936\,b_{{3}}-384\,b_{{1}}
\right)b_{{2}}-4968\,{b_{{3}}}^{2}+\left(5580-984\,b_{{1}}
\right)b_{{3}}\right.\right.\right.\\
\left.\left.\left.+720\,b_{{1}}-48\,{b_{{1}}}^{2}\right) {v}^{2}+1800+
920\,b_{{3}}+320\,b_{{2}}+80\,b_{{1}}\right)\left(\cos\left(v
\right)\right)^{3}\right.\\
\left.+\left(-9936\,\left(b_{{3}}+\frac{2}{23}\,b_{{1}}+{
\frac{8}{23}}\,b_{{2}}\right)\left(\frac{5}{9}+\left( b_{{3}}+\frac{4}{9}\,b_{{2
}}+\frac{1}{9}\,b_{{1}}\right){v}^{2}\right)v\sin\left(v\right)\right.\right.\\ \left.\left.-648\,b
_{{2}}\left(b_{{3}}+\frac{4}{9}\,b_{{2}}+\frac{1}{9}\,b_{{1}}\right) ^{2}{v}^{6}+
9072\,\left({\frac{23}{63}}\,b_{{2}}+b_{{3}}+\frac{1}{9}\,b_{{1}} \right)\right.\right.\\
\left.\left.\left(b_{{3}}+\frac{4}{9}\,b_{{2}}+\frac{1}{9}\,b_{{1}}\right){v}^{4}+\left( -624
\,{b_{{2}}}^{2}+\left(4280-2268\,b_{{3}}-252\,b_{{1}}\right) b_{{2}}\right.\right.\right.\\
\left.\left.\left.-1944\,\left(b_{{3}}+\frac{1}{9}\,b_{{1}}\right)\left( b_{{3}}+\frac{1}{9}\,b_{{
1}}-{\frac{140}{27}}\right)\right) {v}^{2}+2800+260\,b_{{2}}+40\,b
_{{1}}\right.\right.\\
\left.\left.+360\,b_{{3}}\right)\left(\cos\left(v\right)\right)^{2
}+\left(-2592\,\left(b_{{3}}+\frac{1}{9}\,b_{{1}}+{\frac {13}{18}}\,b_{{2}
}\right)\left(\frac{5}{9}\right.\right.\right.\\
\left.\left.\left.+\left(b_{{3}}+\frac{4}{9}\,b_{{2}}+\frac{1}{9}\,b_{{1}}
\right){v}^{2}\right)v\sin\left(v\right)+2187\,\left(b_{{3}}
-\frac{1}{27}\,b_{{1}}\right)\left(b_{{3}}+\frac{4}{9}\,b_{{2}}+\frac{1}{9}\,b_{{1}}
\right)^{2}{v}^{6}\right.\right.\\
\left.\left.-1863\,\left({\frac {212}{207}}\,b_{{2}}+{\frac
{7}{23}}\,b_{{1}}+b_{{3}}\right)\left( b_{{3}}+\frac{4}{9}\,b_{{2}}+\frac{1}{9}\,b_
{{1}}\right){v}^{4}+\left(480\,{b_{{2}}}^{2}\right.\right.\right.\\
\left.\left.\left.+\left( -2120+120\,b_
{{1}}+2520\,b_{{3}}\right)b_{{2}}+3240\,{b_{{3}}}^{2}+\left( -4095+
360\,b_{{1}}\right)b_{{3}}\right.\right.\right.\\
\left.\left.\left.-555\,b_{{1}}\right) {v}^{2}-1325-600\,b_
{{3}}-200\,b_{{2}}\right)\cos\left(v\right)+2160\,\left(\frac{5}{9}+
\left(b_{{3}}+\frac{4}{9}\,b_{{2}}\right.\right.\right.\\
\left.\left.\left.+\frac{1}{9}\,b_{{1}}\right){v}^{2}\right)v
\left(\frac{1}{3}\,b_{{2}}+b_{{3}}\right)\sin\left(v\right) +324\,b_{{2
}}\left(b_{{3}}+\frac{4}{9}\,b_{{2}}+\frac{1}{9}\,b_{{1}}\right) ^{2}{v}^{6}\right.\\
\left.-648\,
\left(-\frac{1}{9}\,b_{{2}}+b_{{3}}+\frac{1}{9}\,b_{{1}}\right)\left( b_{{3}}+\frac{4}{9}
\,b_{{2}}+\frac{1}{9}\,b_{{1}}\right){v}^{4}+\left( 144\,{b_{{2}}}^{2}\right.\right.\\
\left.\left.+\left(-520+228\,b_{{3}}+132\,b_{{1}}\right)b_{{2}}-216\,\left( b_
{{3}}+\frac{1}{9}\,b_{{1}}\right)\left(b_{{3}}-b_{{1}}+{\frac {155}{24}}
\right)\right){v}^{2}\right.\\
\left.-75-60\,b_{{2}}-40\,b_{{1}}+40\,b_{{3}}
\right)\\
/\left(\frac{5}{9}+\left(b_{{3}}+\frac{4}{9}\,b_{{2}}+\frac{1}{9}\,b_{{1}}
\right){v}^{2}\right)^{3}
\ec
The four equations to be solved are:

$$
PL=0, \quad PL'=0, \quad PL''=0, \quad
b_0 = 5-2\,b_{{2}}-2\,b_{{1}}-2\,b_{{3}}
$$

\noindent and the coefficients are given below:

\begin{equation}
\label{coeff_meth_phl_deriv_2}
\begin{array}{l}
\ds b_0=\frac{b_{0,num}}{2 D}, \quad b_1=-\frac{b_{1,num}}{8 D}, \quad b_2=\frac{b_{2,num}}{4 D}, \quad b_3=-\frac{b_{3,num}}{8 D}\vspace{5pt},\\
\mbox{where } D = {{v}^{4}
\left(\sin\left(v\right)\right)^{4}\left(\cos\left(v
\right)-1\right)}\;\;\; \mbox{and}
\end{array}
\end{equation}

\bc
b_{{0,num}}=-6+25\,\left(\cos\left(v\right)\right) ^{3}{v}^{4}+16\,\left(\cos\left(v\right)\right)^{7}{v}^{2}-120\,
\left(\cos\left(v\right)\right)^{4}\\
-32\,\sin\left(v\right)
v\left(\cos\left(v\right)\right)^{6}-96\,\sin\left(v
\right)v\left(\cos\left(v\right)\right)^{7}+32\,\left(\cos
\left(v\right)\right)^{8}{v}^{2}\\
-36\,\cos\left(v\right){v}^{
2}+15\,\cos\left(v\right){v}^{4}+20\,\left(\cos\left(v
\right)\right)^{4}{v}^{4}-96\,\left(\cos\left(v\right)
\right)^{8}\\
+30\,{v}^{4}\left(\cos\left(v\right)\right)^{2}+
20\,\sin\left(v\right)v-12\,{v}^{2}+10\,\left(\cos\left(v
\right)\right)^{5}{v}^{4}\\
+160\,\sin\left(v\right)v\left(\cos
\left(v\right)\right)^{5}+140\,\sin\left(v\right)v\left(
\cos\left(v\right)\right)^{4}\\
-60\,\sin\left(v\right)v
\left(\cos\left(v\right)\right)^{3}-134\,\sin\left(v
\right)v\left(\cos\left(v\right)\right)^{2}+2\,\sin\left(v
\right)v\cos\left(v\right)\\
+18\,\cos\left(v\right)+30\,
\left(\cos\left(v\right)\right)^{2}-54\,\left(\cos\left(v
\right)\right)^{3}+192\,\left(\cos\left(v\right)\right)^{6}\\
+36\,\left(\cos\left(v\right)\right)^{5}+24\,\left(\cos
\left(v\right)\right)^{2}{v}^{2}-64\,\left(\cos\left(v
\right)\right)^{6}{v}^{2}+88\,\left(\cos\left(v\right)
\right)^{3}{v}^{2}\\
-68\,\left(\cos\left(v\right)\right)^{5}{v
}^{2}+20\,\left(\cos\left(v\right)\right)^{4}{v}^{2}
\ec

\bc
b_{{1,num}}=-18-192\,\left(\cos\left(v\right) \right)
^{7}+120\,\left(\cos\left(v\right)\right)^{3}{v}^{4}+128\,
\left(\cos\left(v\right)\right)^{7}{v}^{2}\\
-480\,\left(\cos
\left(v\right)\right)^{4}-320\,\sin\left(v\right)v\left(
\cos\left(v\right)\right)^{6}-192\,\sin\left(v\right)v
\left(\cos\left(v\right)\right)^{7}\\
+64\,\left(\cos\left(v
\right)\right)^{8}{v}^{2}-104\,\cos\left(v\right){v}^{2}+30\,
\cos\left(v\right){v}^{4}+15\,{v}^{4}\\
+60\,\left(\cos\left(v
\right)\right)^{4}{v}^{4}-192\,\left(\cos\left(v\right)
\right)^{8}+75\,{v}^{4}\left(\cos\left(v\right)\right)^{2}+
64\,\sin\left(v\right)v\\
-40\,{v}^{2}+496\,\sin\left(v\right)v
\left(\cos\left(v\right)\right)^{5}+680\,\sin\left(v
\right)v\left(\cos\left(v\right)\right)^{4}\\
-320\,\sin
\left(v\right)v\left(\cos\left(v\right)\right)^{3}-418\,
\sin\left(v\right)v\left(\cos\left(v\right)\right)^{2}+10
\,\sin\left(v\right)v\cos\left(v\right)\\
+42\,\cos\left(v
\right)+162\,\left(\cos\left(v\right)\right)^{2}-258\,
\left(\cos\left(v\right)\right)^{3}+528\,\left(\cos\left(v
\right)\right)^{6}\\
+408\,\left(\cos\left(v\right)\right)^{5
}+32\,\left(\cos\left(v\right)\right)^{2}{v}^{2}-176\,\left(
\cos\left(v\right)\right)^{6}{v}^{2}+336\,\left(\cos\left(v
\right)\right)^{3}{v}^{2}\\
-360\,\left(\cos\left(v\right)
\right)^{5}{v}^{2}+120\,\left(\cos\left(v\right)\right)^{4}{
v}^{2}
\ec

\bc
b_{{2,num}}=-6-96\,\left(\cos\left(v\right)\right) ^{7
}+15\,\left(\cos\left(v\right)\right)^{3}{v}^{4}+48\,\left(
\cos\left(v\right)\right)^{7}{v}^{2}\\
-84\,\left(\cos\left(v
\right)\right)^{4}-128\,\sin\left(v\right)v\left(\cos
\left(v\right)\right)^{6}-40\,\cos\left(v\right){v}^{2}+15\,
\cos\left(v\right){v}^{4}\\
+30\,{v}^{4}\left(\cos\left(v
\right)\right)^{2}+20\,\sin\left(v\right) v-8\,{v}^{2}+48\,\sin
\left(v\right)v\left(\cos\left(v\right)\right)^{5}\\
+240\,
\sin\left(v\right)v\left(\cos\left(v\right)\right)^{4}-48
\,\sin\left(v\right)v\left(\cos\left(v\right)\right)^{3}\\
-126\,\sin\left(v\right)v\left(\cos\left(v\right)\right)^{2
}-6\,\sin\left(v\right)v\cos\left(v\right)+18\,\cos\left(v
\right)+42\,\left(\cos\left(v\right)\right)^{2}\\
-114\,\left(\cos\left(v\right)\right)^{3}+48\,\left(\cos\left(v
\right)\right)^{6}+192\,\left(\cos\left(v\right)\right)^{5
}+16\,\left(\cos\left(v\right)\right)^{2}{v}^{2}\\
+128\,\left(\cos\left(v\right)\right)^{3}{v}^{2}-136\,\left(\cos\left(v
\right)\right)^{5}{v}^{2}-8\,\left(\cos\left(v\right)
\right)^{4}{v}^{2}
\ec

\bc
b_{{3,num}}=48\,\left(\cos\left(v\right) \right)^{6}{
v}^{2}-48\,\left(\cos\left(v\right)\right)^{6}+48\,\left(
\cos\left(v\right)\right)^{5}\\
-80\,\sin\left(v\right)v
\left(\cos\left(v\right)\right)^{5}-48\,\left(\cos\left(v
\right)\right)^{5}{v}^{2}+80\,\sin\left(v\right)v\left(\cos
\left(v\right)\right)^{4}\\
-96\,\left(\cos\left(v\right)
\right)^{4}{v}^{2}+72\,\left(\cos\left(v\right)\right)^{4}+
96\,\left(\cos\left(v\right)\right)^{3}{v}^{2}-78\,\left(
\cos\left(v\right)\right)^{3}\\
+104\,\sin\left(v\right)v
\left(\cos\left(v\right)\right)^{3}-18\,\left(\cos\left(v
\right)\right)^{2}-102\,\sin\left(v\right)v\left(\cos
\left(v\right)\right)^{2}\\
+48\,\left(\cos\left(v\right)
\right)^{2}{v}^{2}+5\,{v}^{4}\left(\cos\left(v\right)\right)
^{2}-18\,\sin\left(v\right)v\cos\left(v\right)-48\,\cos
\left(v\right){v}^{2}\\
+30\,\cos\left(v\right)+10\,\cos\left(v
\right){v}^{4}-6+16\,\sin\left(v\right)v+5\,{v}^{4}
\ec

The Taylor series expansions of the coefficients are given below:

\bc
\ds b_{{0}}=-{\frac {12629}{3024}}+{\frac {45767}{12096}}\,{v}^{2}-{\frac{9837221}{7983360}}\,{v}^{4}+{\frac {153204313}{653837184}}\,{v}^{6}\vspace{5pt}\\
\ds -{\frac {2356782689}{87178291200}}\,{v}^{8}+{\frac {20347993339}{9700566220800}}\,{v}^{10}-{\frac {8744186458121}{77410518441984000}}\,{v}^{12}
+\ldots\\
\\
\ds b_{{1}}={\frac {20483}{4032}}-{\frac {45767}{16128}}\,{v}^{2}+{\frac {2943449}{3548160}}\,{v}^{4}-{\frac {107557349}{792529920}}\,{v}^{6}\vspace{5pt}\\
\ds +{\frac {5074066909}{348713164800}}\,{v}^{8}-{\frac {10190684747}{9484998082560}}\,{v}^{10}+{\frac {5994017812967}{103214024589312000}}\,{v}^{12}-\ldots\\
\\
\ds b_{{2}}=-{\frac {3937}{2016}}+{\frac {45767}{40320}}\,{v}^{2}-{\frac {8607}{39424}}\,{v}^{4}+{\frac {51408821}{2724321600}}\,{v}^{6}\vspace{5pt}\\
\ds -{\frac{35318011}{34871316480}}\,{v}^{8}+{\frac {3348191339}{118562476032000}}\,{v}^{10}-{\frac {56104711163}{43667471941632000}}\,{v}^{12}+\ldots\\
\\
\ds b_{{3}}={\frac {17671}{12096}}-{\frac {45767}{241920}}\,{v}^{2}+{\frac {22153}{4561920}}\,{v}^{4}-{\frac {41092123}{130767436800}}\,{v}^{6}\vspace{5pt}\\
\ds -{\frac {7321421}{348713164800}}\,{v}^{8}-{\frac {5642643317}{2134124568576000}}\,{v}^{10}-{\frac {210863655707}{681212562289459200}}\,{v}^{12}
-
\ec

\subsection{Third optimized method with zero $PL$, $PL'$, $PL''$ and $PL'''$}

All four free coefficients of the third method will be determined by conditions $\{PL=0, PL'=0, PL''=0, PL'''=0\}$.

We use formula \eqref{phl_multi_defn} to compute the phase-lag and then its first, second and third derivative in respect to $v$:

\bc
\ds PL=\left(16\,\left(\cos\left(v\right)\right)^{4}+\left(-16+8
\,b_{{3}}{v}^{2}\right)\left(\cos\left(v\right)\right)^{3}+
\left(4\,b_{{2}}{v}^{2}-8\right)\left(\cos\left(v\right)
\right)^{2}\right.\\
\left.+\left(10+\left(2\,b_{{1}}-6\,b_{{3}}\right){v}^{2}
\right)\cos\left(v\right)-2+\left(-2\,b_{{2}}+b_{{0}}\right)
{v}^{2}\right)/\\
\left(10+\left(18\,b_{{3}}+8\,b_{{2}}+2\,b_{{1}}\right){v}^{2}\right)
\ec

\bc
\ds PL'=\left(-16\,v\left(9\,b_{{3}}+4\,b_{{2}}+b_{{1}}\right)\left(
\cos\left(v\right)\right)^{4}+\left(\left(-160+\left(-128
\,b_{{2}}\right.\right.\right.\right.\\
\left.\left.\left.\left.-32\,b_{{1}}-288\,b_{{3}}\right){v}^{2}\right)\sin
\left(v\right)+16\,\left(4\,b_{{2}}+b_{{1}}+\frac{23}{2}\,b_{{3}}
\right)v\right)\left(\cos\left(v\right)\right)^{3}+\right.\\
\left.\left(-12\,\left(-2+b_{{3}}{v}^{2}\right)\left(5+\left(9\,b_
{{3}}+4\,b_{{2}}+b_{{1}}\right){v}^{2}\right)\sin\left(v
\right)\right.\right.\\
\left.\left.+8\,\left(\frac{13}{2}\,b_{{2}}+9\,b_{{3}}+b_{{1}}\right)v
\right)\left(\cos\left(v\right)\right)^{2}+\left(-4\,
\left(5+\left(9\,b_{{3}}+4\,b_{{2}}+b_{{1}}\right){v}^{2}
\right)\right.\right.\\
\left.\left.\left(-2+b_{{2}}{v}^{2}\right)\sin\left(v\right)-120
\,v\left(b_{{3}}+\frac{1}{3}\,b_{{2}}\right)\right)\cos\left(v
\right)-\left(5+\left(9\,b_{{3}}+4\,b_{{2}}+b_{{1}}\right){v}^
{2}\right)\right.\\
\left.\left(5+\left(b_{{1}}-3\,b_{{3}}\right){v}^{2}
\right)\sin\left(v\right)+5\,v\left(b_{{0}}+{\frac{18}{5}}\,b
_{{3}}-\frac{2}{5}\,b_{{2}}+\frac{2}{5}\,b_{{1}}\right)\right)\\
\left(5+\left(9\,b_{{3}}+4\,b_{{2}}+b_{{1}}\right){v}^{2}\right)^{2}
\ec

\bc
PL'' = {\frac{1}{64}}\,\left(\left(-2048\,\left(\frac{9}{4}\,b_{{3}}+b_{{2}}+1
/4\,b_{{1}}\right)^{2}{v}^{4}+768\,\left(\frac{1}{4}\,b_{{1}}-{\frac{20}{
3}}+b_{{2}}+\frac{9}{4}\,b_{{3}}\right)\right.\right.\\
\left.\left.\left(\frac{9}{4}\,b_{{3}}+b_{{2}}+\frac{1}{4}\,b_{
{1}}\right){v}^{2}-320\,b_{{2}}-80\,b_{{1}}-3200-720\,b_{{3}}
\right)\left(\cos\left(v\right)\right)^{4}\right.\\
\left.+\left(2048\,
\left(\frac{5}{4}+\left(\frac{9}{4}\,b_{{3}}+b_{{2}}+\frac{1}{4}\,b_{{1}}\right){v}^{2}
\right)v\left(\frac{9}{4}\,b_{{3}}+b_{{2}}+\frac{1}{4}\,b_{{1}}\right)\sin
\left(v\right)\right.\right.\\
\left.\left.-576\,b_{{3}}\left(\frac{9}{4}\,b_{{3}}+b_{{2}}+\frac{1}{4}\,b_{{1
}}\right)^{2}{v}^{6}+1152\,\left(\frac{9}{4}\,b_{{3}}+b_{{2}}+\frac{1}{4}\,b_{{1}}
\right)\left(b_{{3}}+b_{{2}}+\frac{1}{4}\,b_{{1}}\right){v}^{4}\right.\right.\\
\left.\left.+\left(
-768\,{b_{{2}}}^{2}+\left(2880-384\,b_{{1}}-3936\,b_{{3}}\right)b_
{{2}}-4968\,{b_{{3}}}^{2}+\left(5580-984\,b_{{1}}\right)b_{{3}}\right.\right.\right.\\
\left.\left.\left.-48
\,{b_{{1}}}^{2}+720\,b_{{1}}\right){v}^{2}+1800+920\,b_{{3}}+80\,b_{
{1}}+320\,b_{{2}}\right)\left(\cos\left(v\right)\right)^{3}\right.\\
+\left.\left(-1536\,\left(\frac{5}{4}+\left(\frac{9}{4}\,b_{{3}}+b_{{2}}+\frac{1}{4}\,b_{{1}}
\right){v}^{2}\right)\left(\frac{1}{4}\,b_{{1}}+{\frac{23}{8}}\,b_{{3}}
+b_{{2}}\right)v\sin\left(v\right)\right.\right.\\
\left.\left.-128\,b_{{2}}\left(\frac{9}{4}\,b_{{
3}}+b_{{2}}+\frac{1}{4}\,b_{{1}}\right)^{2}{v}^{6}+1472\,\left({\frac{63}
{23}}\,b_{{3}}+{\frac{7}{23}}\,b_{{1}}+b_{{2}}\right)\right.\right.\\
\left.\left.\left(\frac{9}{4}\,b
_{{3}}+b_{{2}}+\frac{1}{4}\,b_{{1}}\right){v}^{4}+\left(-624\,{b_{{2}}}^{2
}+\left(-252\,b_{{1}}+4280-2268\,b_{{3}}\right)b_{{2}}\right.\right.\right.\\
\left.\left.\left.-1944\,
\left(b_{{3}}+\frac{1}{9}\,b_{{1}}\right)\left(b_{{3}}+\frac{1}{9}\,b_{{1}}-{
\frac{140}{27}}\right)\right){v}^{2}+2800+360\,b_{{3}}+40\,b_{{1}
}+260\,b_{{2}}\right)\right.\\
\left.\left(\cos\left(v\right)\right)^{2}+
\left(-832\,\left(\frac{5}{4}+\left(\frac{9}{4}\,b_{{3}}+b_{{2}}+\frac{1}{4}\,b_{{1}}
\right){v}^{2}\right)\left(b_{{2}}+{\frac{18}{13}}\,b_{{3}}+\frac{2}{13}\,b_{{1}}\right)v\sin\left(v\right)\right.\right.\\
\left.\left.+432\,\left(b_{{3}}-\frac{1}{27}
\,b_{{1}}\right)\left(\frac{9}{4}\,b_{{3}}+b_{{2}}+\frac{1}{4}\,b_{{1}}\right)^{
2}{v}^{6}-848\,\left({\frac{207}{212}}\,b_{{3}}+b_{{2}}+{\frac{63}
{212}}\,b_{{1}}\right)\right.\right.\\
\left.\left.\left(\frac{9}{4}\,b_{{3}}+b_{{2}}+\frac{1}{4}\,b_{{1}}
\right){v}^{4}+\left(480\,{b_{{2}}}^{2}+\left(-2120+2520\,b_{{3}
}+120\,b_{{1}}\right)b_{{2}}+3240\,{b_{{3}}}^{2}\right.\right.\right.\\
\left.\left.\left.+\left(360\,b_{{1}
}-4095\right)b_{{3}}-555\,b_{{1}}\right){v}^{2}-1325-200\,b_{{2}}-
600\,b_{{3}}\right)\cos\left(v\right)\right.\\
\left.+320\,\left(b_{{2}}+3\,b_
{{3}}\right)\left(\frac{5}{4}+\left(\frac{9}{4}\,b_{{3}}+b_{{2}}+\frac{1}{4}\,b_{{1}}
\right){v}^{2}\right)v\sin\left(v\right)+64\,b_{{2}}\right.\\
\left.\left(\frac{9}{4}\,b_{{3}}+b_{{2}}+\frac{1}{4}\,b_{{1}}\right)^{2}{v}^{6}+32\,\left(\frac{9}{4}\,b
_{{3}}+b_{{2}}+\frac{1}{4}\,b_{{1}}\right)\left(-9\,b_{{3}}+b_{{2}}-b_{{1}
}\right){v}^{4}\right.\\
\left.+\left(24\,{b_{{2}}}^{2}+\left(-220-18\,b_{{1}}-
162\,b_{{3}}-60\,b_{{0}}\right)b_{{2}}-486\,\left(b_{{3}}+\frac{1}{9}\,b_{
{1}}\right)\right.\right.\\
\left.\left.\left(b_{{3}}+{\frac{5}{18}}\,b_{{0}}+{\frac{40}{27}}
+\frac{1}{9}\,b_{{1}}\right)\right){v}^{2}-200+10\,b_{{1}}-10\,b_{{2}}+90
\,b_{{3}}+25\,b_{{0}}\right)\\
/\left(\frac{5}{4}+\left(\frac{9}{4}\,b_{{3}}+b_{{2}}+\frac{1}{4}\,b_{{1}}\right){v}^{2}\right)^{3}
\ec

\bc
PL''' = \left(768\,v\left(\left(9\,b_{{3}}+4\,b_{{2}}+b_{{1}}\right)^{
2}{v}^{4}-\frac{1}{4}\,\left(9\,b_{{3}}+4\,b_{{2}}+b_{{1}}\right)\right.\right.\\
\left.\left.\left(b
_{{1}}-40+9\,b_{{3}}+4\,b_{{2}}\right){v}^{2}+5\,b_{{2}}+\frac{5}{4}\,b_{{1}
}+{\frac{45}{4}}\,b_{{3}}+25\right)\right.\\
\left.\left(9\,b_{{3}}+4\,b_{{2}}+b_
{{1}}\right)\left(\cos\left(v\right)\right)^{4}+\left(512
\,\left(5+\left(9\,b_{{3}}+4\,b_{{2}}+b_{{1}}\right){v}^{2}
\right)\right.\right.\\
\left.\left.\left(\left(9\,b_{{3}}+4\,b_{{2}}+b_{{1}}\right)^{2}{v}
^{4}-{\frac{9}{8}}\,\left(9\,b_{{3}}+4\,b_{{2}}+b_{{1}}\right)
\left(b_{{1}}-{\frac{80}{9}}+9\,b_{{3}}+4\,b_{{2}}\right){v}^{2}\right.\right.\right.\\
\left.\left.\left.+{\frac{135}{8}}\,b_{{3}}+\frac{15}{2}\,b_{{2}}+25+{\frac{15}{8}}\,b_{{1}}
\right)\sin\left(v\right)-432\,\left(\left(9\,b_{{3}}+4\,b_{
{2}}+b_{{1}}\right)^{2}{v}^{4}\right.\right.\right.\\
\left.\left.\left.-\frac{4}{9}\,\left(b_{{1}}-{\frac{45}{2}}+
9\,b_{{3}}+4\,b_{{2}}\right)\left(9\,b_{{3}}+4\,b_{{2}}+b_{{1}}
\right){v}^{2}+{\frac{80}{9}}\,b_{{2}}+{\frac{20}{9}}\,b_{{1}}+20
\,b_{{3}}+25\right)\right.\right.\\
\left.\left.v\left(4\,b_{{2}}+b_{{1}}+\frac{23}{2}\,b_{{3}}
\right)\right)\left(\cos\left(v\right)\right)^{3}+\left(
108\,\left(b_{{3}}\left(9\,b_{{3}}+4\,b_{{2}}+b_{{1}}\right)^{2}
{v}^{6}-2\,\left(4\,b_{{3}}\right.\right.\right.\right.\\
\left.\left.\left.\left.+4\,b_{{2}}+b_{{1}}\right)\left(9\,b_
{{3}}+4\,b_{{2}}+b_{{1}}\right){v}^{4}+\left(414\,{b_{{3}}}^{2}+
\left(328\,b_{{2}}-155+82\,b_{{1}}\right)\right.\right.\right.\right.\\
\left.\left.\left.\left.b_{{3}}+4\,\left(b_{{1}
}+4\,b_{{2}}\right)\left(b_{{1}}-5+4\,b_{{2}}\right)\right){v}
^{2}-50-{\frac{80}{3}}\,b_{{2}}-{\frac{20}{3}}\,b_{{1}}-{\frac{230}
{3}}\,b_{{3}}\right)\right.\right.\\
\left.\left.\left(5+\left(9\,b_{{3}}+4\,b_{{2}}+b_{{1}}
\right){v}^{2}\right)\sin\left(v\right)-672\,v\left(\left(
9\,b_{{3}}+4\,b_{{2}}+b_{{1}}\right)^{2}\right.\right.\right.\\
\left.\left.\left.\left(b_{{1}}+{\frac{61}{
14}}\,b_{{2}}+9\,b_{{3}}\right){v}^{4}-\frac{1}{7}\,\left(9\,b_{{3}}+4\,b_
{{2}}+b_{{1}}\right)\left(81\,{b_{{3}}}^{2}+\right.\right.\right.\right.\\
\left.\left.\left.\left.\left(-630+{\frac{
189}{2}}\,b_{{2}}+18\,b_{{1}}\right)b_{{3}}+26\,{b_{{2}}}^{2}+
\left(\frac{21}{2}\,b_{{1}}-305\right)b_{{2}}+{b_{{1}}}^{2}-70\,b_{{1}}
\right){v}^{2}\right.\right.\right.\\
\left.\left.\left.+{\frac{405}{7}}\,{b_{{3}}}^{2}+\left({\frac{135}{
2}}\,b_{{2}}+225+{\frac{90}{7}}\,b_{{1}}\right)b_{{3}}+{\frac{130}
{7}}\,{b_{{2}}}^{2}+\left(\frac{15}{2}\,b_{{1}}+{\frac{1525}{14}}\right)b
_{{2}}+25\,b_{{1}}\right.\right.\right.\\
\left.\left.\left.+\frac{5}{7}\,{b_{{1}}}^{2}\right)\right)\left(\cos
\left(v\right)\right)^{2}+\left(16\,\left(b_{{2}}\left(9\,
b_{{3}}+4\,b_{{2}}+b_{{1}}\right)^{2}{v}^{6}-14\,\left(9\,b_{{3}}+
4\,b_{{2}}+b_{{1}}\right)\right.\right.\right.\\
\left.\left.\left.\left(9\,b_{{3}}+{\frac{23}{7}}\,b_{{2}}
+b_{{1}}\right){v}^{4}+\left(729\,{b_{{3}}}^{2}+\left(-1260+162
\,b_{{1}}+{\frac{1701}{2}}\,b_{{2}}\right)b_{{3}}+234\,{b_{{2}}}^{2}\right.\right.\right.\right.\\
\left.\left.\left.\left.+\left(-535+{\frac{189}{2}}\,b_{{1}}\right)b_{{2}}-140\,b_{{1}}+
9\,{b_{{1}}}^{2}\right){v}^{2}-350-15\,b_{{1}}-{\frac{195}{2}}\,b_{
{2}}-135\,b_{{3}}\right)\right.\right.\\
\left.\left.\left(5+\left(9\,b_{{3}}+4\,b_{{2}}+b_{{
1}}\right){v}^{2}\right)\sin\left(v\right)+288\,v\left(
\left(b_{{1}}+{\frac{51}{4}}\,b_{{3}}+{\frac{53}{12}}\,b_{{2}}
\right)\right.\right.\right.\\
\left.\left.\left.\left(9\,b_{{3}}+4\,b_{{2}}+b_{{1}}\right)^{2}{v}^{4}-5\,
\left(9\,b_{{3}}+4\,b_{{2}}+b_{{1}}\right)\left(9\,{b_{{3}}}^{2}
+\left(b_{{1}}+7\,b_{{2}}-{\frac{51}{2}}\right)b_{{3}}\right.\right.\right.\right.\\
\left.\left.\left.\left.+\frac{4}{3}\,{b_{{
2}}}^{2}+\left(\frac{1}{3}\,b_{{1}}-{\frac{53}{6}}\right)b_{{2}}-2\,b_{{1
}}\right){v}^{2}+225\,{b_{{3}}}^{2}+\left(25\,b_{{1}}+{\frac{1275
}{4}}+175\,b_{{2}}\right)b_{{3}}\right.\right.\right.\\
\left.\left.\left.+{\frac{100}{3}}\,{b_{{2}}}^{2}+
\left({\frac{25}{3}}\,b_{{1}}+{\frac{1325}{12}}\right)b_{{2}}+25
\,b_{{1}}\right)\right)\cos\left(v\right)+\left(\left(b_{{
1}}-27\,b_{{3}}\right)\right.\right.\\
\left.\left.\left(9\,b_{{3}}+4\,b_{{2}}+b_{{1}}\right)
^{2}{v}^{6}+63\,\left({\frac{23}{7}}\,b_{{3}}+b_{{1}}+{\frac{212}{
63}}\,b_{{2}}\right)\left(9\,b_{{3}}+4\,b_{{2}}+b_{{1}}\right){v
}^{4}\right.\right.\\
\left.\left.+\left(-9720\,{b_{{3}}}^{2}+\left(4095-7560\,b_{{2}}-1080\,b_
{{1}}\right)b_{{3}}-1440\,{b_{{2}}}^{2}\right.\right.\right.\\
\left.\left.\left.+\left(-360\,b_{{1}}+2120
\right)b_{{2}}+555\,b_{{1}}\right){v}^{2}+1800\,b_{{3}}+600\,b_{{2
}}+1325\right)\right.\\
\left.\left(5+\left(9\,b_{{3}}+4\,b_{{2}}+b_{{1}}
\right){v}^{2}\right)\sin\left(v\right)+48\,v\left(\left(
\frac{13}{2}\,b_{{2}}+9\,b_{{3}}+b_{{1}}\right)\right.\right.\\
\left.\left.\left(9\,b_{{3}}+4\,b_{{2}}
+b_{{1}}\right)^{2}{v}^{4}+\frac{1}{2}\,\left(9\,b_{{3}}+4\,b_{{2}}+b_{{1}
}\right)\right.\right.\\
\left.\left.\left(81\,{b_{{3}}}^{2}+\left({\frac{45}{2}}\,b_{{0}}+
180+27\,b_{{2}}+18\,b_{{1}}\right)b_{{3}}-4\,{b_{{2}}}^{2}+\left(3
\,b_{{1}}+130+10\,b_{{0}}\right)b_{{2}}\right.\right.\right.\\
\left.\left.\left.+\left(b_{{1}}+\frac{5}{2}\,b_{{0}}
+20\right)b_{{1}}\right){v}^{2}-{\frac{405}{2}}\,{b_{{3}}}^{2}+
\left(-{\frac{135}{2}}\,b_{{2}}-{\frac{225}{4}}\,b_{{0}}+225-45\,b
_{{1}}\right)b_{{3}}\right.\right.\\
\left.\left.+10\,{b_{{2}}}^{2}+\left(-25\,b_{{0}}+{\frac{
325}{2}}-\frac{15}{2}\,b_{{1}}\right)b_{{2}}-\frac{5}{2}\,\left(-10+\frac{5}{2}\,b_{{0}}+b
_{{1}}\right)b_{{1}}\right)\right)\\
/\left(5+\left(9\,b_{{3}}+4\,b_{{2}}+b_{{1}}\right){v}^{2}\right)^{4}
\ec

After solving the system:

$$
PL=0, \quad PL'=0, \quad PL''=0, \quad PL'''=0
$$

\noindent we get the coefficients:

\begin{equation}
\label{coeff_meth_phl_deriv_3}
\begin{array}{l}
\ds b_0=-\frac{b_{0,num}}{3 D}, \quad b_1=\frac{b_{1,num}}{4 D}, \quad b_2=-\frac{b_{2,num}}{2 D}, \quad b_3=\frac{b_{3,num}}{12 D}\vspace{5pt},\\
\mbox{where } D={{v}^{5}\left(\cos\left(v\right)+1\right)\left(\sin\left(v\right)\right)^{3}}
\;\;\; \mbox{and}
\end{array}
\end{equation}

\bc
b_{{0,num}}=192\,(\cos(v))^{6}
{v}^{2}-126\,\sin(v)(\cos(v)
)^{3}{v}^{3}+99\,\sin(v)(\cos(v
))^{2}v\\
-126\,\sin(v)(\cos
(v))^{2}{v}^{3}-18\,\sin(v){v}^{
3}\cos(v)+630\,\sin(v)(\cos
(v))^{3}v\\
-144\,\sin(v)(
\cos(v))^{6}v+48\,\sin(v){v}^{3}
(\cos(v))^{7}-288\,\sin(v
)(\cos(v))^{7}v\\
-144\,(
\cos(v))^{2}-12\,(\cos(v)
)^{3}+336\,(\cos(v))^{4}+48\,
\sin(v){v}^{3}(\cos(v))^{
6}\\
+30\,{v}^{2}+36\,\cos(v)+144\,(\cos(v
))^{7}{v}^{2}-24\,(\cos(v)
)^{5}+249\,(\cos(v))^{2}{v}^{2}\\
-418\,(\cos(v))^{3}{v}^{2}-662\,(
\cos(v))^{4}{v}^{2}+148\,(\cos(v
))^{5}{v}^{2}-99\,\cos(v)\sin(v
)v\\
-9\,\sin(v)v+66\,\cos(v){v}^{2
}+96\,\sin(v){v}^{3}(\cos(v)
)^{5}+96\,\sin(v)(\cos(v)
)^{4}{v}^{3}\\
-126\,\sin(v)(\cos(v
))^{4}v-18\,\sin(v){v}^{3}-192\,
(\cos(v))^{8}+176\,{v}^{2}(\cos
(v))^{8}\\
-288\,\sin(v)(
\cos(v))^{5}v
\ec

\bc
b_{{1,num}}=12+352\,(\cos(v))^{
6}{v}^{2}-36\,\sin(v)(\cos(v)
)^{3}{v}^{3}+168\,\sin(v)(\cos(v
))^{2}v\\
-100\,\sin(v)(\cos
(v))^{2}{v}^{3}-92\,\sin(v){v}^{
3}\cos(v)+96\,\sin(v)(\cos
(v))^{3}v\\
-672\,\sin(v)(
\cos(v))^{6}v-384\,(\cos(v
))^{7}+36\,(\cos(v))^{2}\\
-48\,(\cos(v))^{3}-48\,(\cos
(v))^{4}+128\,\sin(v){v}^{3}
(\cos(v))^{6}+33\,{v}^{2}-48\,\cos
(v)\\
+448\,(\cos(v))^{7}{v}
^{2}+480\,(\cos(v))^{5}-129\,(
\cos(v))^{2}{v}^{2}-196\,(\cos(v
))^{3}{v}^{2}\\
-316\,(\cos(v)
)^{4}{v}^{2}-464\,(\cos(v))^{5}{
v}^{2}+12\,\cos(v)\sin(v)v-45\,\sin
(v)v\\
+197\,\cos(v){v}^{2}+128\,\sin
(v){v}^{3}(\cos(v))^{5}-
32\,\sin(v)(\cos(v))^{4}{
v}^{3}\\
+504\,\sin(v)(\cos(v)
)^{4}v+4\,\sin(v){v}^{3}-288\,\sin(v
)(\cos(v))^{5}v
\ec

\bc
b_{{2,num}}=152\,(\cos(v))^{6}
{v}^{2}-96\,(\cos(v))^{6}+48\,\sin
(v){v}^{3}(\cos(v))^{5}\\
-192\,\sin(v)(\cos(v))^{5}
v+104\,(\cos(v))^{5}{v}^{2}-72\,\sin
(v)(\cos(v))^{4}v\\
+48\,\sin(v)(\cos(v))^{4}{v}^{
3}-244\,(\cos(v))^{4}{v}^{2}+144\,
(\cos(v))^{4}\\
+216\,\sin(v)(\cos(v))^{3}v-44\,\sin(
v)(\cos(v))^{3}{v}^{3}-134\,
(\cos(v))^{3}{v}^{2}\\
-12\,(\cos(v))^{3}+39\,\sin(v)(\cos
(v))^{2}v-44\,\sin(v)(
\cos(v))^{2}{v}^{3}\\
-48\,(\cos(v))^{2}+79\,(\cos(v))^{2}
{v}^{2}-33\,\cos(v)\sin(v)v-4\,\sin
(v){v}^{3}\cos(v)\\
+18\,\cos(v){v}^{2}+12\,\cos(v)-3\,\sin(v)v
-4\,\sin(v){v}^{3}+10\,{v}^{2}
\ec

\bc
b_{{3,num}}=208\,(\cos(v))^{5}
{v}^{2}-96\,(\cos(v))^{5}-216\,\sin
(v)(\cos(v))^{4}v\\
+120\,(\cos(v))^{4}{v}^{2}+96\,\sin(v
)(\cos(v))^{4}{v}^{3}-360\,
(\cos(v))^{3}{v}^{2}\\
+144\,(\cos
(v))^{3}-72\,\sin(v)(\cos
(v))^{3}v+72\,\sin(v)(
\cos(v))^{3}{v}^{3}\\
+252\,\sin(v)
(\cos(v))^{2}v-157\,(\cos(
v))^{2}{v}^{2}-120\,\sin(v)(\cos
(v))^{2}{v}^{3}\\
-12\,(\cos(v))^{2}+149\,\cos(v){v}^{2}-48\,\cos
(v)+36\,\cos(v)\sin(v)v\\
-72\,\sin(v){v}^{3}\cos(v)-45\,\sin
(v)v+25\,{v}^{2}+24\,\sin(v){v}^{3}+12
\ec

The Taylor series expansions of the coefficients are given below:

\bc
\ds b_{{0}}=-{\frac {12629}{3024}}+{\frac {45767}{9072}}\,{v}^{2}-{\frac {27865393}{11975040}}\,{v}^{4}+{\frac {557684327}{817296480}}\,{v}^{6}\vspace{5pt}\\
\ds -{\frac {235111157089}{1569209241600}}\,{v}^{8}+{\frac {575696865983}{26676557107200}}\,{v}^{10}-{\frac {73845973877087}{32750603956224000}}\,{v}^{12}+...\\
\\
\ds b_{{1}}={\frac {20483}{4032}}-{\frac {45767}{12096}}\,{v}^{2}+{\frac {3549253}{2280960}}\,{v}^{4}-{\frac {36881797}{99066240}}\,{v}^{6}\vspace{5pt}\\
\ds +{\frac {95714204623}{2092278988800}}\,{v}^{8}-{\frac {138581370311}{35568742809600}}\,{v}^{10}+{\frac {106905916402097}{567677135241216000}}\,{v}^{12}\\
\\
\ds b_{{2}}=-{\frac {3937}{2016}}+{\frac {45767}{30240}}\,{v}^{2}-{\frac {3156581}{7983360}}\,{v}^{4}+{\frac {21796097}{681080400}}\,{v}^{6}\vspace{5pt}\\
\ds -{\frac {2365857293}{1046139494400}}\,{v}^{8}-{\frac {102137141}{17784371404800}}\,{v}^{10}-{\frac {3198002983423}{283838567620608000}}\,{v}^{12}\\
\\
\ds b_{{3}}={\frac {17671}{12096}}-{\frac {45767}{181440}}\,{v}^{2}+{\frac {135959}{47900160}}\,{v}^{4}-{\frac {14453093}{16345929600}}\,{v}^{6}\vspace{5pt}\\
\ds -{\frac {90901339}{896690995200}}\,{v}^{8}-{\frac {1564247467}{106706228428800}}\,{v}^{10}-{\frac {3513993676211}{1703031405723648000}}\,{v}^{12}
\ec

It is noteworthy that the Taylor series expansions of all four optimized methods coincide in the constant term and the coefficient of $v^2$ and differ on the coefficients of $v^4$ and for higher powers.

\subsection{Error analysis}

We present the principal term of the local truncation error of the five methods:

Classical method:
\bc
\ds PLTE_{Classical}={\frac{45767}{725760}}\,y^{(10)}{h}^{10}
\ec

Phase fitted method:
\bc
\ds PLTE_{Phase-Fitted}={\frac{45767}{725760}}\left(y^{(10)}+y^{(8)}{\omega}^{2}\right){h}^{10}
\ec

Zero $PL$ and $PL'$ method:
\bc
\ds PLTE_{1st\;deriv}={\frac{45767}{725760}}\left(y^{(10)}+{\omega}^{4}y^{(6)}+2\,{\omega}^{2}y^{(8)}\right){h}^{10}
\ec

Zero $PL$, $PL'$ and $PL''$ method:
\bc
\ds PLTE_{2nd\;deriv}={\frac{45767}{725760}}\left(y^{(10)}+3\,{\omega}^{4}y^{(6)}+3\,{\omega}^{2}y^{(8)}+{\omega}^{6}y^{(4)}\right){h}^{10}
\ec

Zero $PL$, $PL'$, $PL''$ and $PL'''$ method:
\bc
\ds PLTE_{3rd\;deriv}={\frac{45767}{725760}}\left(6\,y^{(6)}{\omega}^{4}+y^{(10)}+4\,y^{(4)}{\omega}^{6}+4\,{\omega}^{2}y^{(8)}+{\omega}^{8}y^{(2)}\right){h}^{10}
\ec

\noindent where $\omega$ is the dominant frequency of the problem. We also present the principal term of the local truncation error of the above methods for the case of the one-dimensional time-independent Schr\"odinger equation:

Classical method:
\bc
\ds PLTE_{Classical}={\frac{1}{725760}}\,{h}^{10}\Big[-45767\,y{E}^{5}+228835\,y\,{E}^{4}\\
+((-2288350\,W''-457670\,(W)^{2})y-915340\,(W')y'){E}^{3}\\
+((457670\,(W)^{3}+3935962\,W^{(4)}+6865050\,W\,''+4576700\,(W')^{2})y\\
+3661360\,(W^{(3)})y'+2746020\,(W')y'W){E}^{2}+((-228835\,(W)^{4}\\
-6865050\,(W)^{2}W''+(-7871924\,W^{(4)}-9153400\,(W')^{2})W-1327243\,W^{(6)}\\
-9656837\,(W'')^{2}-15469246\,(W')W^{(3)})y-14645440\,(W')(W'')y'\\
-7322720\,W(W^{(3)})y'-2837554\,(W^{(5)})y'-2746020\,(W)^{2}(W')y')E\\
+(45767\,(W)^{5}+2288350\,(W)^{3}W''+(3935962\,W^{(4)}+4576700\,(W')^{2})(W)^{2}\\
+(1327243\,W^{(6)}+9656837\,(W'')^{2}+15469246\,(W')W^{(3)})W\\
+2929088\,(W')W^{(5)}+45767\,W^{(8)}+4485166\,(W'')W^{(4)}\\
+10617944\,(W')^{2}W''+2562952\,(W^{(3)})^{2})y+915340\,(W)^{3}(W')y'\\
+3661360\,(W)^{2}(W^{(3)})y'+(2837554\,(W^{(5)})y'+14645440\,(W')(W'')y')W\\
+3661360\,(W')^{3}y'+366136\,(W^{(7)})y'+12814760\,(W'')(W^{(3)})y'\\
+8238060\,(W')(W^{(4)})y'\Big]
\ec

Phase fitted method:
\bc
\ds PLTE_{Phase-Fitted}={\frac{1}{725760}}\,{h}^{10}\Big[(-45767\,{\overline{W}}+45767\,W)y{E}^{4}\\
+((-1281476\,W''+183068\,W{\overline{W}}-183068\,(W)^{2})y-366136\,W'y'){E}^{3}\\
+(((4851302\,W-1006874\,{\overline{W}})W''-274602\,(W)^{2}{\overline{W}}+274602\,(W)^{3}\\
+3203690\,W^{(4)}+3295224\,(W')^{2})y+((-549204\,{\overline{W}}+1647612\,W)W'\\
+2562952\,W^{(3)})y'){E}^{2}+((-8970332\,(W'')^{2}+(2013748\,W{\overline{W}}\\
-5858176\,(W)^{2})W''+(-7871924\,W+1281476\,{\overline{W}})(W')^{2}\\
-14279304\,W'W^{(3)}+183068\,(W)^{3}{\overline{W}}-183068\,(W)^{4}+732272\,W^{(4)}{\overline{W}}\\
-7139652\,WW^{(4)}-1281476\,W^{(6)})y-12448624\,W'W''y'\\
+((1098408\,W{\overline{W}}-2196816\,(W)^{2})W'+1098408\,W^{(3)}{\overline{W}}\\
-6224312\,WW^{(3)}-2562952\,W^{(5)})y')E+((9656837\,W\\
-686505\,{\overline{W}})(W'')^{2}+(-1006874\,(W)^{2}{\overline{W}}+4485166\,W^{(4)}\\
+2288350\,(W)^{3}+10617944\,(W')^{2})W''+(4576700\,(W)^{2}\\
-1281476\,W{\overline{W}})(W')^{2}+(2929088\,W^{(5)}-1189942\,W^{(3)}{\overline{W}}\\
+15469246\,WW^{(3)})W'+45767\,W^{(8)}+45767\,(W)^{5}-45767\,(W)^{4}{\overline{W}}\\
+3935962\,(W)^{2}W^{(4)}+(1327243\,W^{(6)}-732272\,W^{(4)}{\overline{W}})W\\
+2562952\,(W^{(3)})^{2}-45767\,{\overline{W}}\,W^{(6)})y+((-2196816\,{\overline{W}}\\
+14645440\,W)W'+12814760\,W^{(3)})y'W''+(3661360\,(W')^{3}\\
+(8238060\,W^{(4)}+915340\,(W)^{3}-549204\,(W)^{2}{\overline{W}})W'\\
+3661360\,(W)^{2}W^{(3)}+(2837554\,W^{(5)}-1098408\,W^{(3)}{\overline{W}})W\\
-274602\,{\overline{W}}\,W^{(5)}+366136\,W^{(7)})y'\Big]
\ec

Zero $PL$ and $PL'$ method:
\bc
\ds PLTE_{1st\;deriv}={\frac{1}{725760}}\,{h}^{10}\Big[((-45767\,(W)^{2}-594971\,W''-45767\,{{\overline{W}}}^{2}\\
+91534\,W{\overline{W}})y-91534\,W'y'){E}^{3}+((137301\,(W)^{3}-274602\,(W)^{2}{\overline{W}}\\
+(3157923\,W''+137301\,{{\overline{W}}}^{2})W+2517185\,W^{(4)}-1373010\,W''{\overline{W}}\\
+2196816\,(W')^{2})y+(1647612\,W^{(3)}-549204\,W'{\overline{W}}+823806\,W\,W')y'){E}^{2}\\
+((-1235709\,W^{(6)}-137301\,(W)^{4}+274602\,(W)^{3}{\overline{W}}+(-4851302\,W''\\
-137301\,{{\overline{W}}}^{2})(W)^{2}+(-6590448\,(W')^{2}+3386758\,W''{\overline{W}}\\
-6407380\,W^{(4)})W-320369\,W''{{\overline{W}}}^{2}-8283827\,(W'')^{2}+1373010\,{\overline{W}}\,W^{(4)}\\
+2196816\,(W')^{2}{\overline{W}}-13089362\,W'W^{(3)})y+(-2288350\,W^{(5)}\\
-1647612\,(W)^{2}W'+(-5125904\,W^{(3)}+1647612\,W'{\overline{W}})W-274602\,W'{{\overline{W}}}^{2}\\
-10251808\,W'W''+1830680\,{\overline{W}}\,W^{(3)})y')E+(2929088\,W'W^{(5)}\\
+45767\,W^{(8)}+(1327243\,W-91534\,{\overline{W}})W^{(6)}+45767\,(W)^{5}\\
-91534\,(W)^{4}{\overline{W}}+(45767\,{{\overline{W}}}^{2}+2288350\,W'')(W)^{3}+(4576700\,(W')^{2}\\
-2013748\,W''{\overline{W}}+3935962\,W^{(4)})(W)^{2}+(-2562952\,(W')^{2}{\overline{W}}\\
+9656837\,(W'')^{2}+320369\,W''{{\overline{W}}}^{2}-1464544\,{\overline{W}}\,W^{(4)}\\
+15469246\,W'W^{(3)})W+(45767\,{{\overline{W}}}^{2}+4485166\,W'')W^{(4)}\\
+10617944\,(W')^{2}W''-1373010\,(W'')^{2}{\overline{W}}+183068\,(W')^{2}{{\overline{W}}}^{2}\\
+2562952\,(W^{(3)})^{2}-2379884\,W'W^{(3)}{\overline{W}})y+((-549204\,{\overline{W}}\\
+2837554\,W)W^{(5)}+366136\,W^{(7)}+915340\,W'(W)^{3}+(-1098408\,W'{\overline{W}}\\
+3661360\,W^{(3)})(W)^{2}+(-2196816\,{\overline{W}}\,W^{(3)}+14645440\,W'W''\\
+274602\,W'{{\overline{W}}}^{2})W+8238060\,W'W^{(4)}+(12814760\,W''+183068\,{{\overline{W}}}^{2})W^{(3)}\\
-4393632\,W'W''{\overline{W}}+3661360\,(W')^{3})y'\Big]
\ec

Zero $PL$, $PL'$ and $PL''$ method:
\bc
\ds PLTE_{2nd\;deriv}={\frac{1}{725760}}\,{h}^{10}\Big[-183068\,y\,W''{E}^{3}+((45767\,(W)^{3}\\
-137301\,{\overline{W}}\,(W)^{2}+(1784913\,W''+137301\,{{\overline{W}}}^{2})W-1235709\,{\overline{W}}\,W''\\
-45767\,{{\overline{W}}}^{3}+1281476\,(W')^{2}+1876447\,W^{(4)})y+274602\,W'y'W\\
+(915340\,W^{(3)}-274602\,W'{\overline{W}})y'){E}^{2}+((-91534\,(W)^{4}+274602\,(W)^{3}{\overline{W}}\\
+(-274602\,{{\overline{W}}}^{2}-3844428\,W'')(W)^{2}+(-5675108\,W^{(4)}+91534\,{{\overline{W}}}^{3}\\
-5308972\,(W')^{2}+4119030\,{\overline{W}}\,W'')W-11899420\,W'W^{(3)}-7597322\,(W'')^{2}\\
+2746020\,(W')^{2}{\overline{W}}-823806\,{{\overline{W}}}^{2}W''+1922214\,W^{(4)}{\overline{W}}-1189942\,W^{(6)})y\\
-1098408\,(W)^{2}W'y'+(-4027496\,W^{(3)}+1647612\,W'{\overline{W}})y'W\\
+(2196816\,W^{(3)}{\overline{W}}-8054992\,W'W''-549204\,W'{{\overline{W}}}^{2}-2013748\,W^{(5)})y')E\\
+(45767\,(W)^{5}-137301\,(W)^{4}{\overline{W}}+(137301\,{{\overline{W}}}^{2}+2288350\,W'')(W)^{3}\\
+(-3020622\,{\overline{W}}\,W''+3935962\,W^{(4)}-45767\,{{\overline{W}}}^{3}+4576700\,(W')^{2})(W)^{2}\\
+(9656837\,(W'')^{2}+1327243\,W^{(6)}-3844428\,(W')^{2}{\overline{W}}+961107\,{{\overline{W}}}^{2}W''\\
-2196816\,W^{(4)}{\overline{W}}+15469246\,W'W^{(3)})W+45767\,W^{(8)}-137301\,W^{(6)}{\overline{W}}\\
+2929088\,W'W^{(5)}+(4485166\,W''+137301\,{{\overline{W}}}^{2})W^{(4)}+2562952\,(W^{(3)})^{2}\\
-3569826\,W'W^{(3)}{\overline{W}}-2059515\,(W'')^{2}{\overline{W}}+(-45767\,{{\overline{W}}}^{3}\\
+10617944\,(W')^{2})W''+549204\,(W')^{2}{{\overline{W}}}^{2})y+915340\,(W)^{3}W'y'\\
+(-1647612\,W'{\overline{W}}+3661360\,W^{(3)})y'(W)^{2}+(14645440\,W'W''\\
-3295224\,W^{(3)}{\overline{W}}+823806\,W'{{\overline{W}}}^{2}+2837554\,W^{(5)})y'W\\
+366136\,W^{(7)}y'+(-823806\,W^{(5)}{\overline{W}}+8238060\,W'W^{(4)}\\
+(549204\,{{\overline{W}}}^{2}+12814760\,W'')W^{(3)}-91534\,W'{{\overline{W}}}^{3}\\
-6590448\,W'W''{\overline{W}}+3661360\,(W')^{3})y'\Big]
\ec

Zero $PL$, $PL'$, $PL''$ and $PL'''$ method:
\bc
\ds PLTE_{3rd\;deriv}={\frac{1}{725760}}\,{h}^{10}\Big[((1281476\,W^{(4)}+(-732272\,{\overline{W}}\\
+732272\,W)W''+549204\,(W')^{2})y+366136\,W^{(3)}y'){E}^{2}+((-1144175\,W^{(6)}\\
+(2379884\,{\overline{W}}-4942836\,W)W^{(4)}-10709478\,W'W^{(3)}-6910817\,(W'')^{2}\\
+(4210564\,W{\overline{W}}-1373010\,{{\overline{W}}}^{2}-2837554\,(W)^{2})W''\\
-45767\,(W)^{4}+183068\,(W)^{3}{\overline{W}}-274602\,(W)^{2}{{\overline{W}}}^{2}\\
+(-4027496\,(W')^{2}+183068\,{{\overline{W}}}^{3})W+2929088\,(W')^{2}{\overline{W}}\\
-45767\,{{\overline{W}}}^{4})y+(-1739146\,W^{(5)}+(2196816\,{\overline{W}}\\
-2929088\,W)W^{(3)}+1098408\,W{\overline{W}}\,W'-5858176\,W'W''-549204\,W'{{\overline{W}}}^{2}\\
-549204\,(W)^{2}W')y')E+(45767\,W^{(8)}+(1327243\,W-183068\,{\overline{W}})W^{(6)}\\
+2929088\,W'W^{(5)}+(3935962\,(W)^{2}+4485166\,W''+274602\,{{\overline{W}}}^{2}\\
-2929088\,W{\overline{W}})W^{(4)}+2562952\,(W^{(3)})^{2}+(-4759768\,W'{\overline{W}}\\
+15469246\,W\,W')W^{(3)}+(-2746020\,{\overline{W}}+9656837\,W)(W'')^{2}\\
+(2288350\,(W)^{3}+1922214\,W{{\overline{W}}}^{2}-4027496\,(W)^{2}{\overline{W}}\\
+10617944\,(W')^{2}-183068\,{{\overline{W}}}^{3})W''+45767\,(W)^{5}\\
-183068\,(W)^{4}{\overline{W}}+274602\,(W)^{3}{{\overline{W}}}^{2}+(-183068\,{{\overline{W}}}^{3}\\
+4576700\,(W')^{2})(W)^{2}+(-5125904\,(W')^{2}{\overline{W}}+45767\,{{\overline{W}}}^{4})W\\
+1098408\,(W')^{2}{{\overline{W}}}^{2})y+(366136\,W^{(7)}+(-1098408\,{\overline{W}}\\
+2837554\,W)W^{(5)}+8238060\,W'W^{(4)}+(3661360\,(W)^{2}\\
-4393632\,W{\overline{W}}+1098408\,{{\overline{W}}}^{2}+12814760\,W'')W^{(3)}+(-8787264\,W'{\overline{W}}\\
+14645440\,W\,W')W''+1647612\,W{{\overline{W}}}^{2}W'-2196816\,(W)^{2}{\overline{W}}\,W'\\
+3661360\,(W')^{3}+915340\,(W)^{3}W'-366136\,{{\overline{W}}}^{3}W')y'\Big]
\ec

The principal terms of the local truncation errors presented above are collected in respect to the energy $E$ in descending order. As we can easily see, the maximum power of $E$ in the error for each case is:

\begin{itemize}
\item $E^5$ for the classical method
\item $E^4$ for the phase-fitted method
\item $E^3$ for the zero $PL$ and $PL'$ method
\item $E^3$ for the zero $PL$, $PL'$ and $PL''$ method and
\item $E^2$ for the zero $PL$, $PL'$, $PL''$ and $PL'''$ method.
\end{itemize}

A low maximum power of $E$ is crucial when integrating the Schr\"odinger equation using a high value of energy.

\subsection{Stability analysis}
The stability analysis of the methods concerns the application of the test problem $y''=-\omega y$.

Here we present the characteristic equations of the five methods:

\bc
C.E._{Classical}=1+{\lambda}^{8}+{\frac{1}{12096}}\,(17671\,{s}^{2}-24192){\lambda}^{7}+{\frac{1}{12096}}\,(-23622\,{s}^{2}\\
+24192){\lambda}^{6}+{\frac{1}{12096}}\,(61449\,{s}^{2}-12096){\lambda}^{5}-{\frac{12629}{3024}}\,{s}^{2}{\lambda}^{4}\\
+{\frac{1}{12096}}\,(61449\,{s}^{2}-12096){\lambda}^{3}+{\frac{1}{12096}}\,(-23622\,{s}^{2}+24192){\lambda}^{2}\\
+{\frac{1}{12096}}\,(17671\,{s}^{2}-24192)\lambda
\ec

\bc
C.E._{Phase-Fitted}=-{\frac{109}{32}}\,(1+{\lambda}^{2}-2\,\lambda\,\cos(s))(-{\frac{32}{109}}\,(\lambda-1)^{6}(\cos(s))^{3}\\
+({\frac{96}{109}}+{\frac{96}{109}}\,{\lambda}^{6}+({s}^{2}-{\frac{416}{109}}){\lambda}^{5}+(-{\frac{808}{327}}\,{s}^{2}+{\frac{880}{109}}){\lambda}^{4}+({\frac{1202}{327}}\,{s}^{2}-{\frac{1120}{109}}){\lambda}^{3}\\
+(-{\frac{808}{327}}\,{s}^{2}+{\frac{880}{109}}){\lambda}^{2}+({s}^{2}-{\frac{416}{109}})\lambda)(\cos(s))^{2}+(-{\frac{96}{109}}-{\frac{96}{109}}\,{\lambda}^{6}\\
+(-{\frac{404}{327}}\,{s}^{2}+{\frac{296}{109}}){\lambda}^{5}+(-{\frac{480}{109}}+{\frac{736}{327}}\,{s}^{2}){\lambda}^{4}+({\frac{560}{109}}-{\frac{1144}{327}}\,{s}^{2}){\lambda}^{3}\\
+(-{\frac{480}{109}}+{\frac{736}{327}}\,{s}^{2}){\lambda}^{2}+(-{\frac{404}{327}}\,{s}^{2}+{\frac{296}{109}})\lambda)\cos(s)+{\frac{32}{109}}+{\frac{32}{109}}\,{\lambda}^{6}\\
+(-{\frac{72}{109}}+{\frac{137}{327}}\,{s}^{2}){\lambda}^{5}+({\frac{80}{109}}-{\frac{56}{109}}\,{s}^{2}){\lambda}^{4}+(-{\frac{80}{109}}+{\frac{302}{327}}\,{s}^{2}){\lambda}^{3}\\
+({\frac{80}{109}}-{\frac{56}{109}}\,{s}^{2}){\lambda}^{2}+(-{\frac{72}{109}}+{\frac{137}{327}}\,{s}^{2})\lambda)(\cos(s)-1)^{-3}
\ec

\bc
C.E._{1st\,Deriv}={\frac{125}{48}}\,(-{\frac{96}{125}}\,\lambda\,s(\lambda-1)^{4}(\cos(s))^{5}+{\frac{48}{125}}\,(4\,\lambda\,\sin(s)\\
+s({\lambda}^{2}+4\,\lambda+1))(\lambda-1)^{4}(\cos(s))^{4}+(-{\frac{192}{125}}\,\lambda\,(\lambda-1)^{4}\sin(s)\\
-2\,({\frac{48}{125}}+{\frac{48}{125}}\,{\lambda}^{6}-{\frac{192}{125}}\,{\lambda}^{5}+({\frac{396}{125}}+{s}^{2}){\lambda}^{4}+(-{\frac{38}{25}}\,{s}^{2}-{\frac{504}{125}}){\lambda}^{3}\\
+({\frac{396}{125}}+{s}^{2}){\lambda}^{2}-{\frac{192}{125}}\,\lambda)s)(\cos(s))^{3}+\lambda\,(-{\frac{96}{125}}\,(\lambda-1)^{4}\sin(s)+(({s}^{2}-{\frac{132}{125}}){\lambda}^{4}\\
+(-{\frac{62}{25}}\,{s}^{2}+{\frac{648}{125}}){\lambda}^{3}+({\frac{98}{25}}\,{s}^{2}-{\frac{1032}{125}}){\lambda}^{2}+(-{\frac{62}{25}}\,{s}^{2}+{\frac{648}{125}})\lambda+{s}^{2}\\
-{\frac{132}{125}})s)(\cos(s))^{2}+({\frac{24}{25}}\,\lambda\,(\lambda-1)^{4}\sin(s)+{\frac{12}{25}}\,s(\frac{8}{5}+\frac{8}{5}\,{\lambda}^{6}+(-{\frac{24}{5}}+{s}^{2}){\lambda}^{5}\\
+(\frac{1}{6}\,{s}^{2}+{\frac{34}{5}}){\lambda}^{4}+(-{\frac{36}{5}}-\frac{1}{3}\,{s}^{2}){\lambda}^{3}+(\frac{1}{6}\,{s}^{2}+{\frac{34}{5}}){\lambda}^{2}+(-{\frac{24}{5}}+{s}^{2})\lambda))\cos(s)\\
-{\frac{24}{125}}\,\lambda\,(\lambda-1)^{4}\sin(s)-{\frac{13}{25}}\,(1+{\lambda}^{2})({\frac{48}{65}}+{\frac{48}{65}}\,{\lambda}^{4}+({s}^{2}-{\frac{132}{65}}){\lambda}^{3}\\
+(-{\frac{14}{13}}\,{s}^{2}+{\frac{168}{65}}){\lambda}^{2}+({s}^{2}-{\frac{132}{65}})\lambda)s)(1+{\lambda}^{2}-2\,\lambda\,\cos(s)){s}^{-1}\\
(\cos(s)+1)^{-1}(\cos(s)-1)^{-3}
\ec

\bc
C.E._{2nd\,Deriv}=-\frac{5}{8}\,(-{\frac{32}{5}}\,{\lambda}^{2}(\lambda-1)^{2}({s}^{2}-3)(\cos(s))^{7}\\
+{\frac{32}{5}}\,({\lambda}^{2}{s}^{2}-\frac{3}{2}\,{\lambda}^{2}+3\,\lambda\,\sin(s)s+\lambda\,{s}^{2}-3\,\lambda+{s}^{2}-\frac{3}{2})\lambda\,(\lambda-1)^{2}(\cos(s))^{6}\\
-\frac{8}{5}\,(10\,s({\lambda}^{2}+1+\frac{6}{5}\,\lambda)\lambda\,\sin(s)+{\lambda}^{4}{s}^{2}+(4\,{s}^{2}-6){\lambda}^{3}+(18-6\,{s}^{2}){\lambda}^{2}\\
+(4\,{s}^{2}-6)\lambda+{s}^{2})(\lambda-1)^{2}(\cos(s))^{5}+(16\,({\lambda}^{2}-\frac{7}{5}\,\lambda+1)s\lambda\,(\lambda-1)^{2}\sin(s)\\
+\frac{8}{5}\,{\lambda}^{6}{s}^{2}+({\frac{72}{5}}-16\,{s}^{2}){\lambda}^{5}+({\frac{88}{5}}\,{s}^{2}+{\frac{12}{5}}){\lambda}^{4}+(-{\frac{168}{5}}-{\frac{32}{5}}\,{s}^{2}+4\,{s}^{4}){\lambda}^{3}\\
+({\frac{88}{5}}\,{s}^{2}+{\frac{12}{5}}){\lambda}^{2}+({\frac{72}{5}}-16\,{s}^{2})\lambda+\frac{8}{5}\,{s}^{2})(\cos(s))^{4}\\
-4\,(-{\frac{26}{5}}\,({\lambda}^{2}+{\frac{25}{26}}\,\lambda+1)s\lambda\,\sin(s)-\frac{4}{5}\,{\lambda}^{4}{s}^{2}+({\frac{39}{10}}-{\frac{16}{5}}\,{s}^{2}){\lambda}^{3}+({s}^{4}-\frac{9}{5}){\lambda}^{2}\\
+({\frac{39}{10}}-{\frac{16}{5}}\,{s}^{2})\lambda-\frac{4}{5}\,{s}^{2})(\lambda-1)^{2}(\cos(s))^{3}+(-{\frac{102}{5}}\,({\lambda}^{2}-\frac{2}{17}\,\lambda+1)s\lambda\,\\
(\lambda-1)^{2}\sin(s)-{\frac{16}{5}}\,{\lambda}^{6}{s}^{2}+({\frac{64}{5}}\,{s}^{2}-{\frac{18}{5}}+{s}^{4}){\lambda}^{5}+(-8\,{s}^{4}-16\,{s}^{2}-{\frac{24}{5}}){\lambda}^{4}\\
+({\frac{84}{5}}+{\frac{64}{5}}\,{s}^{2}+6\,{s}^{4}){\lambda}^{3}+(-8\,{s}^{4}-16\,{s}^{2}-{\frac{24}{5}}){\lambda}^{2}+({\frac{64}{5}}\,{s}^{2}-{\frac{18}{5}}+{s}^{4})\lambda\\
-{\frac{16}{5}}\,{s}^{2})(\cos(s))^{2}+2\,(-\frac{9}{5}\,({\lambda}^{2}-\frac{4}{9}\,\lambda+1)s\lambda\,\sin(s)-\frac{4}{5}\,{\lambda}^{4}{s}^{2}\\
+(-{\frac{16}{5}}\,{s}^{2}+3+{s}^{4}){\lambda}^{3}+(-\frac{8}{5}\,{s}^{2}+\frac{6}{5}){\lambda}^{2}+(-{\frac{16}{5}}\,{s}^{2}+3+{s}^{4})\lambda\\
-\frac{4}{5}\,{s}^{2})(\lambda-1)^{2}\cos(s)+{\frac{16}{5}}\,s({\lambda}^{2}-\frac{1}{2}\,\lambda+1)\lambda\,(\lambda-1)^{2}\sin(s)\\
+(1+{\lambda}^{2})(\frac{8}{5}\,{\lambda}^{4}{s}^{2}+({s}^{4}-\frac{6}{5}-{\frac{16}{5}}\,{s}^{2}){\lambda}^{3}+({\frac{12}{5}}+{\frac{16}{5}}\,{s}^{2}){\lambda}^{2}\\
+({s}^{4}-\frac{6}{5}-{\frac{16}{5}}\,{s}^{2})\lambda+\frac{8}{5}\,{s}^{2}))(1+{\lambda}^{2}-2\,\lambda\,\cos(s)){s}^{-2}\\
(\cos(s)+1)^{-2}(\cos(s)-1)^{-3}
\ec

\bc
C.E._{3rd\,Deriv}=-(1+{\lambda}^{2}-2\,\lambda\,\cos(s))((-{\frac{88}{3}}\,{s}^{2}+32){\lambda}^{3}(\cos(s))^{7}\\
-8\,(\lambda\,s({s}^{2}-6)\sin(s)+(4-{\frac{31}{6}}\,{s}^{2}){\lambda}^{2}+3\,\lambda\,{s}^{2}+4-{\frac{31}{6}}\,{s}^{2}){\lambda}^{2}\\
(\cos(s))^{6}+12\,(s((-5+{s}^{2}){\lambda}^{2}+(2-\frac{2}{3}\,{s}^{2})\lambda-5+{s}^{2})\lambda\,\sin(s)\\
+(-{\frac{13}{9}}\,{s}^{2}+\frac{2}{3}){\lambda}^{4}+\frac{8}{3}\,{\lambda}^{3}{s}^{2}+({\frac{7}{9}}\,{s}^{2}-\frac{8}{3}){\lambda}^{2}+\frac{8}{3}\,\lambda\,{s}^{2}-{\frac{13}{9}}\,{s}^{2}+\frac{2}{3})\lambda\,(\cos(s))^{5}\\
-6\,(s(({s}^{2}-3){\lambda}^{4}+(-2\,{s}^{2}+4){\lambda}^{3}+(\frac{2}{3}\,{s}^{2}+2){\lambda}^{2}+(-2\,{s}^{2}+4)\lambda\\
-3+{s}^{2})\sin(s)+\frac{5}{3}\,{\lambda}^{4}{s}^{2}+(-8+{\frac{31}{3}}\,{s}^{2}){\lambda}^{3}+(-\frac{2}{3}-{\frac{11}{9}}\,{s}^{2}){\lambda}^{2}\\
+(-8+{\frac{31}{3}}\,{s}^{2})\lambda+\frac{5}{3}\,{s}^{2})\lambda\,(\cos(s))^{4}+(s({\lambda}^{6}{s}^{2}+(6-6\,{s}^{2}){\lambda}^{5}+(66-9\,{s}^{2}){\lambda}^{4}\\
+(-4\,{s}^{2}-3){\lambda}^{3}+(66-9\,{s}^{2}){\lambda}^{2}+(6-6\,{s}^{2})\lambda+{s}^{2})\sin(s)+(30\,{s}^{2}-12){\lambda}^{5}\\
+(-{\frac{245}{6}}\,{s}^{2}-4){\lambda}^{4}+(-8+{\frac{145}{3}}\,{s}^{2}){\lambda}^{3}+(-{\frac{245}{6}}\,{s}^{2}-4){\lambda}^{2}+(30\,{s}^{2}-12)\lambda)\\
(\cos(s))^{3}+(s({\lambda}^{6}{s}^{2}+(-21+6\,{s}^{2}){\lambda}^{5}+(-9\,{s}^{2}+{\frac{27}{2}}){\lambda}^{4}+(12\,{s}^{2}-39){\lambda}^{3}\\
+(-9\,{s}^{2}+{\frac{27}{2}}){\lambda}^{2}+(-21+6\,{s}^{2})\lambda+{s}^{2})\sin(s)+({\frac{157}{12}}\,{s}^{2}+1){\lambda}^{5}\\
+({\frac{44}{3}}\,{s}^{2}-16){\lambda}^{4}+({\frac{173}{6}}\,{s}^{2}-2){\lambda}^{3}+({\frac{44}{3}}\,{s}^{2}-16){\lambda}^{2}+({\frac{157}{12}}\,{s}^{2}+1)\lambda)(\cos(s))^{2}\\
+(-s({\lambda}^{6}{s}^{2}+(3-6\,{s}^{2}){\lambda}^{5}+(3\,{s}^{2}+9){\lambda}^{4}+(-12\,{s}^{2}+3){\lambda}^{3}+(3\,{s}^{2}+9){\lambda}^{2}\\
+(3-6\,{s}^{2})\lambda+{s}^{2})\sin(s)+(4-{\frac{149}{12}}\,{s}^{2}){\lambda}^{5}+({\frac{29}{6}}\,{s}^{2}+4){\lambda}^{4}+(8-{\frac{161}{6}}\,{s}^{2}){\lambda}^{3}\\
+({\frac{29}{6}}\,{s}^{2}+4){\lambda}^{2}+(4-{\frac{149}{12}}\,{s}^{2})\lambda)\cos(s)-({\lambda}^{4}{s}^{2}-{\frac{15}{4}}\,{\lambda}^{3}\\
+(\frac{3}{2}+2\,{s}^{2}){\lambda}^{2}-{\frac{15}{4}}\,\lambda+{s}^{2})s(1+{\lambda}^{2})\sin(s)+(-1-{\frac{25}{12}}\,{s}^{2}){\lambda}^{5}+5\,{\lambda}^{4}{s}^{2}\\
+(-2-{\frac{37}{6}}\,{s}^{2}){\lambda}^{3}+5\,{\lambda}^{2}{s}^{2}+(-1-{\frac{25}{12}}\,{s}^{2})\lambda){s}^{-3}(\cos(s)+1)^{-1}(\sin(s))^{-3}
\ec

From the characteristic equations we evaluate $s_0$ and the interval of periodicity $[0,s_0^2]$. These are given below:

\begin{itemize}
\item $s_0=0.754$ ($[0,0.569]$) for the classical method
\item $s_0=0.803$ ($[0,0.645]$) for the phase-fitted method
\item $s_0=0.874$ ($[0,0.763]$) for the zero $PL$ and $PL'$ method
\item $s_0=1.010$ ($[0,1.020]$) for the zero $PL$, $PL'$ and $PL''$ method and
\item $s_0=1.865$ ($[0,3.478]$) for the zero $PL$, $PL'$, $PL''$ and $PL'''$ method.
\end{itemize}

As we can see, by requiring higher derivatives of the phase-lag to be vanished, we increase the interval of periodicity, which is a very important property.

\section{Numerical results}
\label{Numerical_results}

\subsection{The problems}

The efficiency of the two newly constructed methods will be measured through the integration of two real initial value
problems with oscillating solutions.

\subsubsection{The Schr\"odinger equation}
\label{Intro}

The radial Schr\"{o}dinger equation is given by:

\begin{equation}
\label{Schrodinger}
    y''(x) = \left( \frac{l(l+1)}{x^{2}}+V(x)-E \right) y(x)
\end{equation}

\noindent where $\frac{l(l+1)}{x^{2}}$ is the \textit{centrifugal potential}, $V(x)$ is the \textit{potential}, $E$ is
the \textit{energy} and $W(x) = \frac{l(l+1)}{x^{2}} + V(x)$ is the \textit{effective potential}. It is valid that
${\mathop {\lim} \limits_{x \to \infty}} V(x) = 0$ and therefore ${\mathop {\lim} \limits_{x \to \infty}} W(x) = 0$.

We consider $E>0$ and divide $[0,\infty)$ into subintervals $[a_{i},b_{i}]$ so that $W(x)$ is a constant with value
${\mathop{W_{i}}\limits^{\_}}$. After this the problem \eqref{Schrodinger} can be expressed by the approximation

\begin{equation}
\begin{array}{l}
\label{Schrodinger_simplified}
y''_{i} = ({\mathop{W}\limits^{\_}} - E)\,y_{i}, \quad\quad \mbox{whose solution is}\\
y_{i}(x) = A_{i}\,\exp{\left(\sqrt{{\mathop{W}\limits^{\_}}-E}\,x\right)} +
B_{i}\,\exp{\left(-\sqrt{{\mathop{W}\limits^{\_}}-E}\,x\right)}, \\
 A_{i},\,B_{i}\,\in {\mathbb R}.
\end{array}
\end{equation}

We will integrate problem \eqref{Schrodinger} with $l=0$ at the interval $[0,15]$ using the well known Woods-Saxon
potential

\begin{eqnarray}
\label{Woods_Saxon} V(x) = \frac{u_{0}}{1+q} + \frac{u_{1}\,q}{(1+q)^2}, \quad\quad q =
\exp{\left(\frac{x-x_{0}}{a}\right)}, \quad
\mbox{where}\\
\nonumber u_{0}=-50, \quad a=0.6, \quad x_{0}=7 \quad \mbox{and} \quad u_{1}=-\frac{u_{0}}{a}
\end{eqnarray}

\noindent and with boundary condition $y(0)=0$.

\noindent The potential $V(x)$ decays more quickly than $\frac{l\,(l+1)}{x^2}$, so for large $x$ (asymptotic region) the Schr\"{o}dinger equation \eqref{Schrodinger} becomes

\begin{equation}
\label{Schrodinger_reduced}
    y''(x) = \left( \frac{l(l+1)}{x^{2}}-E \right) y(x)
\end{equation}

\noindent The last equation has two linearly independent solutions $k\,x\,j_{l}(k\,x)$ and\\
$k\,x\,n_{l}(k\,x)$, where $j_{l}$ and $n_{l}$ are the \textit{spherical Bessel} and \textit{Neumann} functions. When $x
\rightarrow \infty$ the solution takes the asymptotic form

\begin{equation}
\label{asymptotic_solution}
\begin{array}{l}
 y(x) \approx A\,k\,x\,j_{l}(k\,x) - B\,k\,x\,n_{l}(k\,x) \\
\approx D[sin(k\,x - \pi\,l/2) + \tan(\delta_{l})\,\cos{(k\,x - \pi\,l/2)}],
\end{array}
\end{equation}

\noindent where $\delta_{l}$ is called \textit{scattering phase shift} and it is given by the following expression:

\begin{equation}
\tan{(\delta_{l})} = \frac{y(x_{i})\,S(x_{i+1}) - y(x_{i+1})\,S(x_{i})} {y(x_{i+1})\,C(x_{i}) - y(x_{i})\,C(x_{i+1})},
\end{equation}

\noindent where $S(x)=k\,x\,j_{l}(k\,x)$, $C(x)=k\,x\,n_{l}(k\,x)$ and $x_{i}<x_{i+1}$ and both belong to the asymptotic
region. Given the energy we approximate the phase shift, the accurate value of which is $\pi/2$ for the above problem.

We will use three different values for the energy:

\begin{itemize}
\item $E_1=989.701916$
\item $E_2=341.495874$
\item $E_3=163.215341$
\end{itemize}

As for the frequency $w$ we will use the suggestion of Ixaru and Rizea \cite{ix_ri}:

\begin{equation}
\omega =
        \begin{cases}
        \sqrt{E-50},  &   x\in[0,\,6.5]\\
        \sqrt{E},     &   x\in[6.5,\,15]
        \end{cases}
\end{equation}

\subsubsection{The N-Body Problem}
The N-body problem is the problem that concerns the movement of N bodies under Newton's law of gravity. It is expressed by a system of vector differential equations

\begin{equation}\label{IVP_N_body}
\begin{array}{l} \displaystyle
\overrightarrow{\ddot{y}_i} = G\,\sum\limits_{j=1,\,j\neq0}^{N}{\frac{m_j\,(\overrightarrow{y_j}-\overrightarrow{y_i})}{|\overrightarrow{y_j}-\overrightarrow{y_i}|^3}}, \quad i=1,2,..,N
\end{array}
\end{equation}

\noindent where $G$ is the gravitational constant, $m_j$ is the mass of body $j$ and $\overrightarrow{y_i}$ is the vector of the position of body $i$.

It is easy to see that each vector differential equation of (\ref{IVP_N_body}) can be analyzed into three simplified differential equations, that express the three directions $x, y, z$. So $\overrightarrow{y_j}-\overrightarrow{y_i}$ expresses the difference between the coordinates of bodies $j$ and $i$ for the corresponding direction, while $|\overrightarrow{y_j}-\overrightarrow{y_i}|$ represents the distance between bodies $i$ and $j$.

The above system of ODEs cannot be solved analytically. Instead we produce a highly accurate numerical solution by using a 10-stage implicit Runge-Kutta method of Gauss with 20th algebraic order, that is also symplectic and A-stable. The method can be easily reproduced using simplifying assumptions for the order conditions (see \cite{butcher}).

The reference solution is obtained by using the previous method to integrate the N-body problem for a specific time-span and for different step-lengths.

In order to find the step-length $h_{opt}$ that gives the best approximation, we have to keep in mind that the total error of a numerical method that integrates a system of ODEs consists of the error due to the truncation error of the method and the roundoff error of all computations. While the global truncation error of the method tends to zero, while $h$ decreases, the opposite happens to the roundoff, which tends to infinity.

If $y_{acc}$ is the analytical solution for a specific time-span of the problem, then let $\epsilon_n=||\overline{y}_{h_n}-y_{acc}||$ and $\varepsilon_n=||\overline{y}_{h_{n+1}}-\overline{y}_{h_{n}}||$, where $\overline{y}_{h_n}$ is the approximate solution of $y$ using a step-length $h_n$. $\epsilon_n$ represents the actual error of the approximation and $\varepsilon_n$ is the best known approximation to the actual error, being the difference of two approximations with different step-lengths. We see that, when $h_n \rightarrow h_{opt}$ $\Rightarrow$ $\epsilon_n \rightarrow \epsilon_{min}$ and $\varepsilon_n \rightarrow \varepsilon_{min}$. The minimum values of the errors $\epsilon_{min}$ and $\varepsilon_{min}$ are positive numbers and depend on the software that is used for the integration and the computer system that it runs on. We can also see that $\epsilon_n$ and $\varepsilon_n$ have similar behavior around $n_{opt}$, meaning that they increase and decrease simultaneously. According to these we find the step-length $h_{opt}$ that minimizes $\varepsilon_n$, which is easily calculated for every $h_n$.

In \cite{hairer} the data for the five outer planet problem is given. This system consists of the sun and the five most distant planets of the solar system. In Table \ref{table_five_outer_planets} we can see the masses, the initial position components and the initial velocity components of the six bodies. Masses are relative to the sun, so that the sun has mass 1. In the computations the sun with the four inner planets are considered one body, so the mass is larger than one. Distances are in astronomical units, time is in earth days and the gravitational constant is $G=2.95912208286 \cdot 10^{-4}$.

\begin{equation}
\begin{array}{c|c|c|c}
\label{table_five_outer_planets}
Planet	&	Mass	&	Initial\; Position	&	Initial\; Velocity	\\
\hline							
Sun	&	1.00000597682	&	0	&	0	\\
	&		            &	0	&	0	\\
	&		            &	0	&	0	\\
\hline							
Jupiter	&	0.000954786104043	&	-3.5023653	&	 ~~0.00565429	\\
	&		                    &	-3.8169847	&	 -0.00412490	\\
	&		                    &	-1.5507963	&	 -0.00190589	\\
\hline							
Saturn	&	0.000285583733151	&	~~9.0755314	&	 ~~0.00168318	\\
	&		                    &	-3.0458353	&	 ~~0.00483525	\\
	&		                    &	-1.6483708	&	 ~~0.00192462	\\
\hline							
Uranus	&	0.0000437273164546	&	~~8.3101420	&	 ~~0.00354178	\\
	&		                    &	-16.2901086	&	 ~~0.00137102	\\
	&		                    &	-7.2521278	&	 ~~0.00055029	\\
\hline							
Neptune	&	0.0000517759138449	&	~~11.4707666	&	 ~~0.00288930	\\
	&		                    &	-25.7294829	&	 ~~0.00114527	\\
	&		                    &	-10.8169456	&	 ~~0.00039677	\\
\hline							
Pluto	&	1/(1.3\cdot10^8)	&	-15.5387357	&	 ~~0.00276725	\\
	&		                &	-25.2225594	&	-0.00170702	 \\
	&		                &	-3.1902382	&	-0.00136504	 \\
\end{array}
\end{equation}

The system of equations (\ref{IVP_N_body}) has been solved for $t\in[0,10^6]$, for which time-span, the previously mentioned method of Gauss produces a $10.5$ decimal digits solution.

We have used $\omega=0.00145044732989$, which is the dominant frequency of the problem, as evaluated by the square root of the spectral radius of matrix A, if the problem is expressed in the form $y''=Ay+B$.

\subsection{The methods}

\begin{itemize}
\item The classical method developed by Quinlan and Tremaine \cite{qt8}
\item The phase-fitted method developed by Panopoulos, Anastassi and\\ Simos \cite{panopoulos_match}
\item The zero $PL$ and $PL'$ method developed here
\item The zero $PL$, $PL'$ and $PL''$ developed here
\item The zero $PL$, $PL'$, $PL''$ and $PL'''$ developed here
\end{itemize}

\subsection{Comparison}

\begin{figure}[htbp]
    \includegraphics[width=\textwidth]{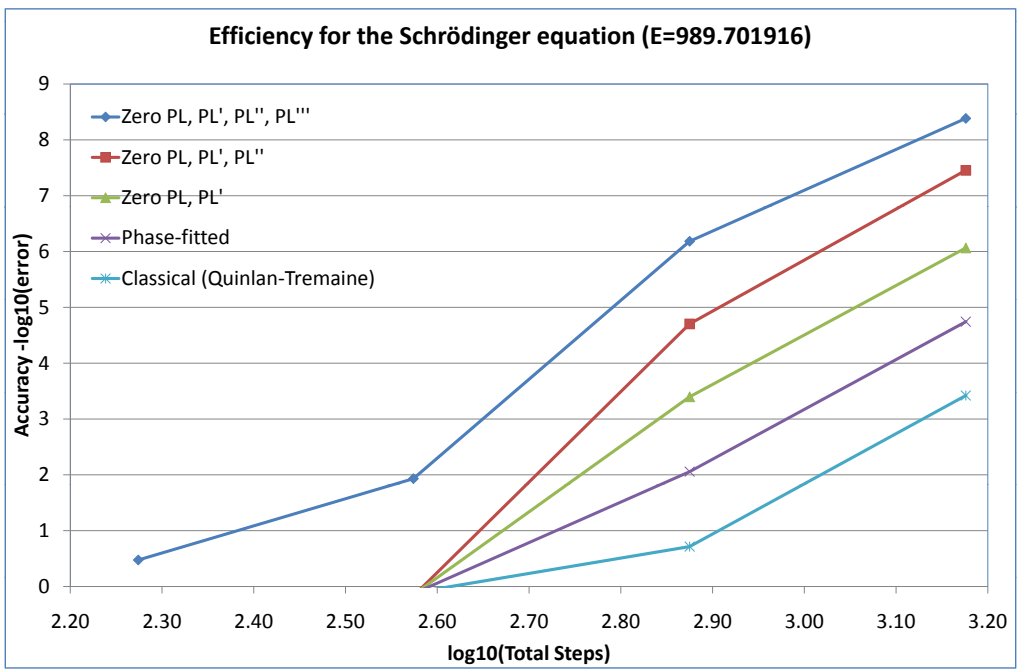}
    \caption{Efficiency for the Schr\"odinger equation using E = 989.701916}
    \label{fig_res_989}
\end{figure}

\begin{figure}[htbp]
    \includegraphics[width=\textwidth]{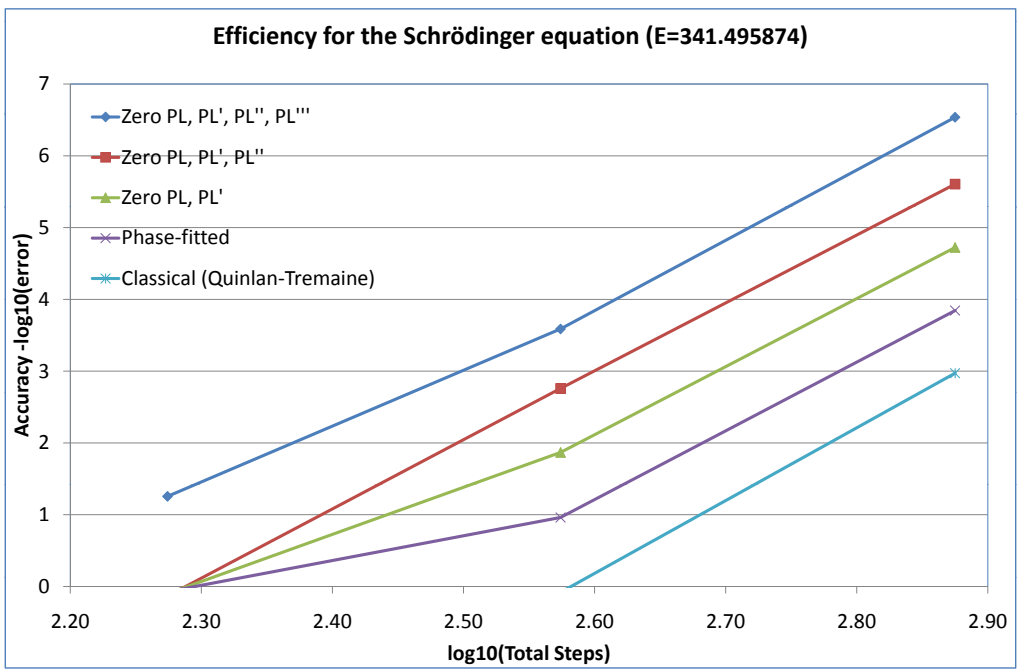}
    \caption{Efficiency for the Schr\"odinger equation using E = 341.495874}
    \label{fig_res_341}
\end{figure}

\begin{figure}[htbp]
    \includegraphics[width=\textwidth]{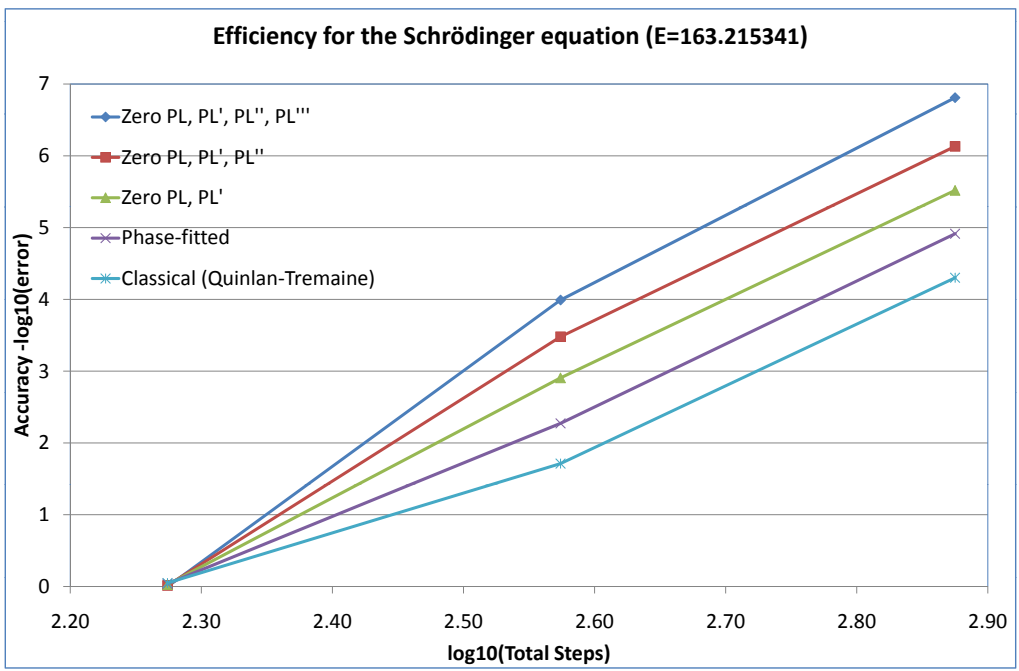}
    \caption{Efficiency for the Schr\"odinger equation using E = 163.215341}
    \label{fig_res_163}
\end{figure}

\begin{figure}[htbp]
    \includegraphics[width=\textwidth]{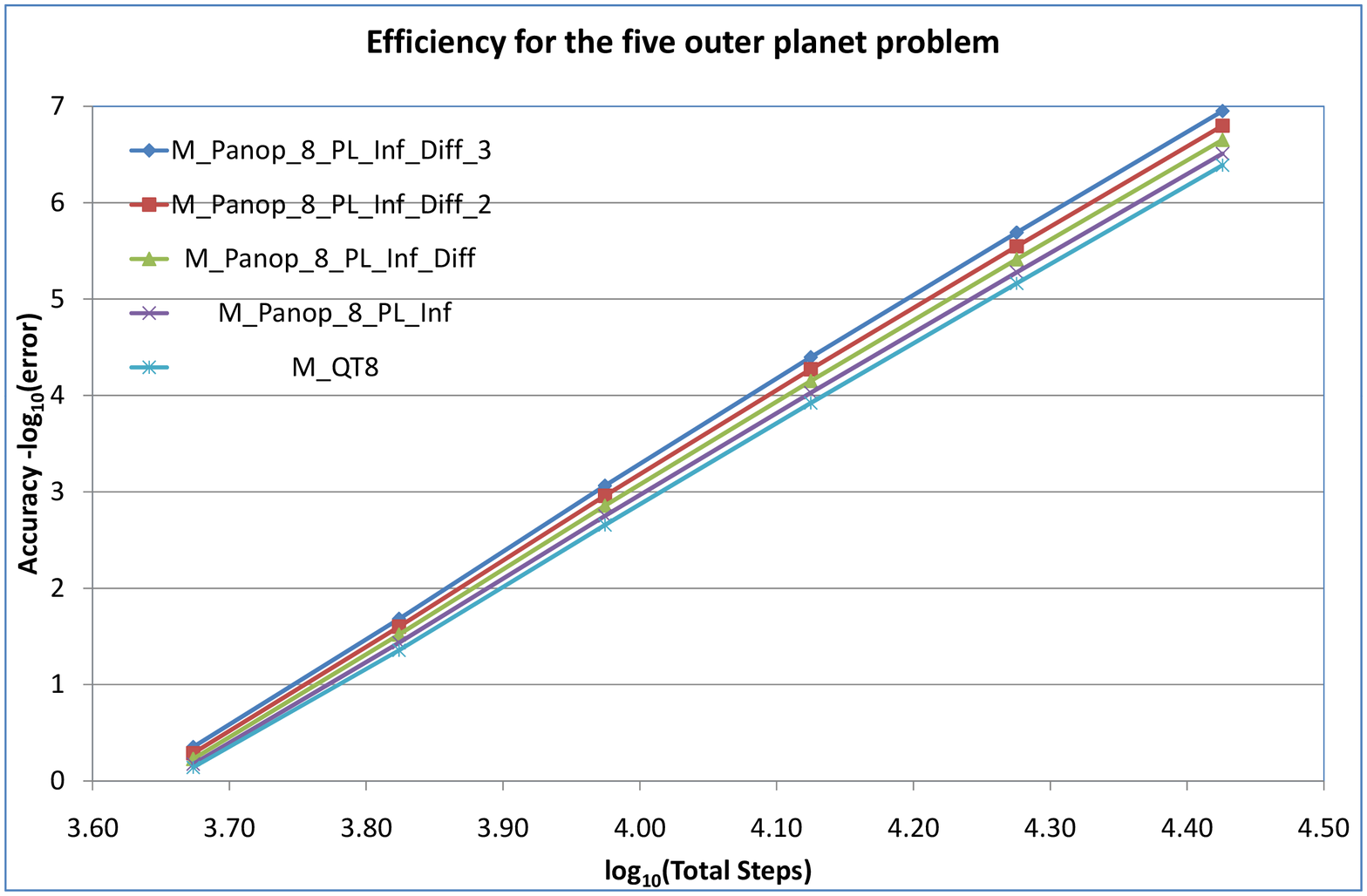}
    \caption{Efficiency for the N-body problem}
    \label{fig_nbody}
\end{figure}

We are presenting the accuracy of the methods expressed by $-\log_{10}$(error at the end point) versus the $\log_{10}$(total steps). In Figures \ref{fig_res_989}, \ref{fig_res_341} and \ref{fig_res_163} we are presenting the efficiency of the methods for the Schr\"odinger equation using a value for the energy equal to i) $989.701916$, ii) $341.495874$ and iii) $163.215341$. Also in Figure \ref{fig_nbody} we present the efficiency for the N-body problem and particularly the five outer planet problem.

We see that for each successive derivative of the phase-lag nullified, we gain in efficiency for both IVPs tested here.

\section{Conclusions}
We have developed three new optimized eight-step symmetric methods with zero phase-lag and derivatives. We showed that the more derivatives of the phase-lag are vanished, the bigger the interval of periodicity and the higher the efficiency of the method. This is the case for both problems tested here. Also the local error truncation analysis shows the relation of the error to the energy, revealing the importance of nullified phase-lag derivatives when integrating the Schr\"odinger equation, especially when using high value of energy.

\end{document}